\documentclass[11pt]{article}

\usepackage{amssymb,latexsym}
\usepackage{amsmath}
\usepackage{citesort}

\title{Locality in GNS Representations of Deformation Quantization}

\author{{\bf
          Stefan
          Waldmann\thanks{Stefan.Waldmann@physik.uni-freiburg.de}
         } \\[1cm]
         Fakult\"{a}t f\"{u}r Physik \\
         Universit\"{a}t Freiburg \\
         Hermann-Herder-Str. 3 \\
         79104 Freiburg i.~Br., F.~R.~G
       }

\date{FR-THEP-99/2 \\[1mm]
      March 1999}

\newcommand{\LS} [1] {{#1} (\!(\lambda)\!)}
\newcommand{\NP} [1] {{#1} \langle\!\langle \lambda^* \rangle\!\rangle}
\newcommand{\CNP} [1] {{#1} \langle\!\langle \lambda \rangle\!\rangle}

\newcommand{\im} {{\mathrm i}}

\newcommand{\Lie} {\mathcal L}
\newcommand{\cc} [1] {\overline {{#1}}}
\newcommand{\supp} {\mathop{{\mathrm {supp}}}}
\newcommand{\lsupp} {{\mathrm {supp}}_\lambda}
\newcommand{\id} {{\mathsf {id}}}
\newcommand{\tr} {{\mathsf {tr}}}
\newcommand{\ad} {{\mathrm {ad}}}
\newcommand{\Hom} {{\mathsf {Hom}}}
\newcommand{\End} {{\mathsf {End}}}
\newcommand{\field} [1] {{\mathsf {{#1}}}}
\newcommand{\SP} [1] {{\left\langle {{#1}} \right\rangle}}

\newcommand{\Exp} {\mathop{{\mathrm {Exp}}}}

\newcommand{\Diff} {{\mathfrak D}}
\newcommand{\Loc} {{\mathfrak L}}
\newcommand{\Bounded} {{\mathfrak B}}
\newcommand{\Unitary} {{\mathfrak U}}
\newcommand{\BLoc} {{\mathfrak {LB}}}

\newcommand{\wick} {{\mathop {*_{\mbox{\rm \tiny Wick}}}}}
\newcommand{\std}  {{\mathop {*_{\mbox{\rm \tiny Std}}}}}
\newcommand{\weyl} {{\mathop {*_{\mbox{\rm \tiny Weyl}}}}}
\newcommand{\wrep} {{\mathop {\varrho_{\mbox{\rm \tiny Weyl}}}}}

\newcommand{\KMS} {{\mbox{\tiny\rm KMS}}}

\newcommand{\AL} {\mathcal A_{\mathsf L}}
\newcommand{\AR} {\mathcal A_{\mathsf R}}
\newcommand{\ALR} {\mathcal A_{\mathsf {LR}}}

\newtheorem{lemma} {Lemma} [section]
\newtheorem{proposition} [lemma] {Proposition}
\newtheorem{theorem} [lemma] {Theorem}
\newtheorem{corollary} [lemma] {Corollary}
\newtheorem{definition}[lemma] {Definition}
\newtheorem{example}[lemma] {Example}
\newtheorem{remark}[lemma]{Remark}

\newenvironment{proof}{\small{\sc Proof:}}{{\hspace*{\fill} $\square$\\}}
\newenvironment{innerproof}{\small {\sc Proof:}}{{\hspace*{\fill}
                            $\bigtriangledown$}}

\numberwithin{equation}{section}

\begin{document}

\maketitle

\begin{abstract}
In the framework of deformation quantization we apply the formal GNS
construction to find representations of the deformed algebras in
pre-Hilbert spaces over $\mathbb C[[\lambda]]$ and establish the
notion of local operators in these pre-Hilbert spaces. The commutant
within the local operators is used to distinguish `thermal' from
`pure' representations. The computation of the local commutant is
exemplified in various situations leading to the physically reasonable
distinction between thermal representations and pure ones. Moreover,
an analogue of von Neumann's double commutant theorem is proved in the
particular situation of a GNS representation with respect to a KMS 
functional and for the Schr\"odinger representation on cotangent
bundles. Finally we prove a formal version of the Tomita-Takesaki
theorem.
\end{abstract}

\newpage

\tableofcontents

\newpage

%
%

\section{Introduction}
\label{IntroSec}

The concept of deformation quantization and star products has been
introduced by Bayen, Flato, Fr{\o}nsdal, Lichnerowicz, and Sternheimer
in \cite{BFFLS78} and is now a well-established and successful
quantization procedure. The central object of deformation quantization
is a star product $*$, a formal associative deformation of the classical
Poisson algebra of complex-valued functions $C^\infty (M)$
on a symplectic or, more
generally, a Poisson manifold $M$, such that in the first order of the
formal parameter $\lambda$ the commutator of the star product yields
$\im$ times the Poisson bracket. Hence $\lambda$ is to be identified
with Planck's constant $\hbar$ whenever the formal series
converges, and the star product algebra $(C^\infty (M)[[\lambda]], *)$
is viewed as observable
algebra of the quantized system corresponding to the classical system
described by the Poisson manifold. 
The existence of such deformations has been shown by
DeWilde and Lecomte \cite{DL83b}, Fedosov \cite{Fed94a,Fed96}, Omori,
Maeda, and Yoshioka \cite{OMY91} in the symplectic case and recently
by Kontsevich \cite{Kon97b} in the general case of a Poisson
manifold. Moreover, the star products have been classified up to
equivalence by Nest and Tsygan \cite{NT95a,NT95b}, Bertelson, Cahen, and
Gutt \cite{BCG97}, Fedosov \cite{Fed96}, Weinstein and Xu 
\cite{WX97}, and Kontsevich \cite{Kon97b}.

Since the algebra structure is by now quite well-understood the
question is raised, how one can encode the notion of states of this
algebra. In order to find a physically reasonable notion of states 
we have introduced in \cite{BW98a} together
with Bordemann a formal analogue of the positive functionals 
and their GNS representation as analogue of the 
well-known construction from $C^*$-algebra theory, see e.g.\ 
\cite{BR87,Con94,Haa93}. In various examples this approach has been
successfully applied
\cite{BW98a,BNW98a,BNW99a,BNPW98,BRW98a,BRW98b}.

While the GNS construction in principle works in a quite general
framework (for $^*$-algebras over ordered rings) we shall discuss in
this article aspects of these representations which are particular to
deformation quantization. Here the observable algebra has an
additional structure since the star products are \emph{local}
in the sense that any bilinear operator of the formal star product
series is a local or even bidifferential operator. This leads to a
`net structure' of the observable algebra similar to the net structure
of the observable algebras in algebraic quantum field theory
\cite{Haa93}. It turns out that the spaces $C^\infty_0 (O)[[\lambda]]$
are two-sided ideals for any open subset $O$ of the manifold $M$. Thus
we want to transfer this `locality structure' to the GNS
representation space and use it to study the GNS representations more
closely.

The main results of this work are organized as follows: 
After a brief summary of deformation quantization in
Section~\ref{BasSec}, the 
first crucial observation in Section~\ref{ComSec} is that one can
assign to any vector in the 
GNS pre-Hilbert space a \emph{support} on $M$. Hence one can think of
these abstract equivalence classes as located on the underlying
manifold. In particular the GNS pre-Hilbert space inherits a net
structure of orthogonal subspaces indexed by the open sets of $M$
where two such subspaces are orthogonal if the corresponding open sets
on the 
manifold are disjoint. This allows for the definition of 
\emph{local operators} where an endomorphism of the GNS representation
space is called local if it is compatible with the above net
structure. It turns out that the GNS representation itself always
yields such local operators. Thus we consider the 
\emph{commutant in the local operators} of the GNS representation in
order to study the 
question when a representation is to be regarded as a `thermal' one or
a `pure' one (Section~\ref{ComSec}).

In Section~\ref{FaithSec} we consider faithful positive functionals
and their GNS 
representations. It turns out that a positive functional is faithful
if and only if its support is the whole manifold and a GNS
representation is faithful if and only if the corresponding positive 
functional is faithful, a feature which is quite different from usual
$C^*$-algebra theory.

In Section~\ref{ExSec} we discuss particular examples of positive
functionals as 
traces and KMS functionals on connected symplectic manifolds, 
$\delta$-functionals on K\"ahler manifolds,
and Schr\"odinger functionals on cotangent bundles. It turns out that
the abstract notion of support of vectors in the GNS pre-Hilbert space
coincides with the usual notion of support in that cases where the GNS
pre-Hilbert space is isomorphic to spaces of formal wave functions. 
Furthermore in these examples the local commutant yields the
physically correct characterization of `thermal' vs.\ `pure'
representations.

In Section~\ref{NeuSec} we define analogues of strong and $^*$-strong
operator 
topologies for the local operators and arrive at an analogue of von
Neumann's double commutant theorem for the cases of GNS
representations with respect to a KMS functional and a Schr\"odinger
functional, respectively. 
Here the completion in the strong operator topology
coincides with the double commutant within the local operators.

Finally, in Section~\ref{ToTaSec} we prove by a simple algebraic
computation a 
formal analogue of the Tomita-Takesaki theorem for the local operators
in a GNS representation of a KMS functional.

In Appendix~\ref{GNSApp} and \ref{FormalApp} we briefly summarize
well-known results on the 
formal GNS construction as well as on formal series and their
$\lambda$-adic topology.

%
%

\section{Basic definitions}
\label{BasSec}

In this section we recall some basic features of deformation
quantization to set-up our notation and discuss the definition of the
support of linear functionals.

Throughout this article $M$ denotes a symplectic or Poisson manifold
endowed with a local or even differential star product $*$. Hereby a
star product is a formal associative deformation of the pointwise
multiplication of the smooth complex-valued functions $C^\infty (M)$
on $M$ in direction of the Poisson bracket. More precisely, 
$*: C^\infty (M)[[\lambda]] \times C^\infty (M) [[\lambda]] 
\to C^\infty (M)[[\lambda]]$ 
is an associative $\mathbb C[[\lambda]]$-bilinear product such that
for $f,g \in C^\infty (M)$
\begin{equation}
\label{StarprodDef}
    f * g = \sum_{r=0}^\infty \lambda^r M_r (f, g),
\end{equation}
with local or even bidifferential operators $M_r$ fulfilling 
$M_0 (f, g) = fg$ and 
$M_1 (f, g) - M_1 (g, f) = \im \{f, g\}$. Moreover we require that for
$r\ge 1$ the operator $M_r$ vanishes on constants, whence 
$1*f = f = f*1$. Note that the $\mathbb C[[\lambda]]$-bilinearity
implies the form (\ref{StarprodDef}) and in particular the
$\lambda$-adically continuity of $*$, see
e.g.\ \cite[Prop.~2.1]{DL88}. With this normalization the formal
parameter $\lambda$ is directly to be identified with Planck's
constant $\hbar$ and may be substituted in convergent
situations. Since we are interested in GNS representations we need a
$^*$-involution, i.e.\ a $\mathbb C[[\lambda]]$-anti-linear involutive
anti-automorphism of the star product $*$. Thus we additionally
require the property 
\begin{equation}
\label{CCinvolution}
    \cc {f * g} = \cc g * \cc f,
\end{equation}
where $f \mapsto \cc f$ denotes the pointwise complex conjugation and
$\lambda$  is considered to be real, i.e.\ we define 
$\cc \lambda := \lambda$. Note that such star products always exist.

Next we consider certain sub-algebras of $C^\infty (M)[[\lambda]]$
indexed by the open subsets of $M$. Viewing elements of 
$C^\infty (M)[[\lambda]]$ as $\mathbb C[[\lambda]]$-valued functions
we define the \emph{support} $\supp f$ of 
$f = \sum_{r=0}^\infty \lambda^r f_r \in C^\infty (M)[[\lambda]]$ 
in the usual way to be the closure of the set 
$\{x \in M \;|\; f(x) \ne 0\}$ which coincides with the 
closure of the set $\bigcup_{r=0}^\infty \supp f_r$.
Note that if $\supp f$ is compact then
$\supp f_r$ is compact for all $r$ but the converse is not true in
general. Due to the required locality of the star product we have
\begin{equation}
\label{localProduct}
    \supp (f*g) \subseteq \supp f \cap \supp g
\end{equation}
for all $f, g \in C^\infty (M)[[\lambda]]$. Furthermore, for any open
set $O \subseteq M$ the spaces 
$\mathcal A (O) = \{f \in C^\infty(M)[[\lambda]] \; | \; 
\supp f \subseteq O \}$ 
and $C^\infty_0 (O)[[\lambda]]$ are two-sided ideals of 
$C^\infty (M)[[\lambda]]$ which are stable under complex conjugation. 
Clearly, if $O \subseteq O'$ then 
$\mathcal A(O) \subseteq \mathcal A(O')$ and
$C^\infty_0 (O)[[\lambda]] \subseteq C^\infty_0 (O')[[\lambda]]$,
respectively. Note that $\mathcal A(O)$ and 
$C^\infty_0 (O)[[\lambda]]$ in general have no unit
element. Note furthermore that if $f \in C^\infty_0 (O)[[\lambda]]$ 
it may happen
that $\supp f$ is not contained in $O$ but only in the closure of $O$.
Moreover, $\supp f$ needs not to be compact. This may cause some
subtleties later on and hence we shall define 
\begin{equation}
\label{NiceGuysDef}
    \mathcal A_0 (O) := 
    \{ f \in C^\infty_0 (O)[[\lambda]] 
    \; | \; \supp f \subseteq O \}
    \subseteq C^\infty_0 (M)[[\lambda]],
\end{equation}
which is again a two-sided ideal in the whole algebra. All these net
structures motivate to speak of a \emph{net of local observables}
similar to algebraic quantum field theory, see e.g.\ \cite{Haa93}. But
note that here locality means locality in phase space and not in
Minkowski space.

In a next step we consider $\mathbb C[[\lambda]]$-linear functionals 
$\omega: C^\infty_0(M)[[\lambda]] \to \mathbb C[[\lambda]]$ or
$\omega: C^\infty (M)[[\lambda]] \to \mathbb C[[\lambda]]$,
respectively. Since in the most
relevant examples we have in mind the functionals incorporate
integrations over $M$ it is reasonable to consider mainly functionals
defined on $C^\infty_0 (M)[[\lambda]]$ only. We define the support 
$\supp \omega$ of such a functional as usual to be the complement of
the union of those open sets $O \subseteq M$ with 
$\omega|_{C^\infty_0 (O)[[\lambda]]}= 0$. Since any 
$\mathbb C[[\lambda]]$-linear functional is automatically continuous
in the $\lambda$-adic topology we obviously have 
$\omega|_{C^\infty_0 (O)[[\lambda]]} = 0$ if and only if
$\omega|_{C^\infty_0 (O)} = 0$. Hence it suffices to `test' $\omega$
on $C^\infty_0 (O)$ in order to determine the support of $\omega$ and
thus we have
\begin{equation}
\label{FunctSuppDef}
    \supp \omega \; := \; 
    M \setminus \mathop{\mbox{$\bigcup$}}
    \limits_{\omega|_{C^\infty_0(O)[[\lambda]]} = 0} O
    \; = \; 
    M \setminus \mathop{\mbox{$\bigcup$}}
    \limits_{\omega|_{C^\infty_0 (O)} = 0} O,
\end{equation}
where $O$ ranges over the open subsets of $M$. Then the following
lemma is obtained completely analogously to the non-formal case of
distributions, see e.g.\ \cite[p.~164]{Rud91}, and does not yet use the
star product. 
\begin{lemma}
\label{BasicSuppLem}
Let $M$ be a manifold and let
$\omega, \omega' : C^\infty_0 (M)[[\lambda]] \to \mathbb C[[\lambda]]$
be $\mathbb C[[\lambda]]$-linear functionals. Then $\supp \omega$ is
closed and 
$\omega|_{C^\infty_0 (M \setminus \supp \omega)[[\lambda]]} = 0$. 
If $f \in C^\infty_0 (M)[[\lambda]]$ fulfills 
$\supp f \cap \supp \omega = \emptyset$ then $\omega (f) = 0$. Finally 
$\supp (\omega + \omega') \subseteq \supp \omega \cup \supp \omega'$
and $\supp \alpha \omega \subseteq \supp \omega$ for 
$\alpha \in \mathbb C[[\lambda]]$.
\end{lemma}
Similar to the case of distributions with compact support we can
construct an extension of a $\mathbb C[[\lambda]]$-linear functional 
$\omega: C^\infty_0 (M)[[\lambda]] \to \mathbb C[[\lambda]]$ to a 
$\mathbb C[[\lambda]]$-linear functional 
$\hat\omega: C^\infty (M)[[\lambda]] \to \mathbb C[[\lambda]]$
provided the support of $\omega$ is compact.
\begin{proposition}
Let $M$ be a manifold and let
$\omega: C^\infty_0 (M)[[\lambda]] \to \mathbb C[[\lambda]]$ be a 
$\mathbb C[[\lambda]]$-linear functional with compact support. Then
there exists a unique extension 
$\hat\omega: C^\infty (M)[[\lambda]] \to \mathbb C[[\lambda]]$ of
$\omega$ with the property that for $f \in C^\infty (M)[[\lambda]]$
with $\supp f \cap \supp \omega = \emptyset$ one has 
$\hat \omega (f) = 0$. 
\end{proposition}
\begin{proof}
Since $\supp \omega$ is assumed to be compact we find a smooth
partition of unity $\chi_0 + \chi_1 = 1$ such that 
$\supp \chi_0 \subseteq M\setminus \supp \omega$,
$\chi_1$ has compact support, and $\chi_1$ is equal to $1$ in an open
neighborhood of $\supp \omega$. Then one verifies easily that 
$\hat\omega (f) := \omega (\chi_1 f)$ is a well-defined extension
having the desired properties which proves the existence. Now let
$\hat\omega'$ be another such extension then for all 
$f \in C^\infty (M)[[\lambda]]$ one has  
$\supp (\chi_0 f) \cap \supp \omega = \emptyset$. Hence
$\hat\omega' (f) = \hat\omega' (\chi_0 f + \chi_1 f) =
\hat\omega'(\chi_1 f) = \omega (f) = \hat \omega (f)$
since $\chi_1 f$ has compact support and clearly
$\omega (\chi_1 f) = \omega (f)$.
\end{proof}

%
%

\section{Locality in GNS representations and commutants}
\label{ComSec}

Let us now consider \emph{positive} $\mathbb C[[\lambda]]$-linear
functionals and their induced GNS representations and investigate the
relations between the support of these functionals and their GNS
representations. For a detailed exposition concerning the GNS
construction in deformation quantization see \cite{BW98a,BNW99a}, and
see App.~\ref{GNSApp} for a short summary and notation.

Let $\omega: C^\infty_0 (M)[[\lambda]] \to \mathbb C[[\lambda]]$ be a
positive $\mathbb C[[\lambda]]$-linear functional. Then
$\mathcal J_\omega$ denotes the Gel'fand ideal of $\omega$ and
$\mathfrak H_\omega$ the 
GNS pre-Hilbert space. Note that the GNS representation 
$\pi_\omega: C^\infty_0 (M)[[\lambda]] 
\to \Bounded (\mathfrak H_\omega)$ extends 
to the whole algebra $C^\infty (M)[[\lambda]]$ since 
$C^\infty_0 (M)[[\lambda]]$ is a two-sided ideal in 
$C^\infty (M)[[\lambda]]$ stable under complex conjugation,  
see \cite[Cor.~1]{BW98a}.

The crucial observation is that we can associate to any vector
$\psi_f$ in the GNS pre-Hilbert space $\mathfrak H_\omega$ a `support'
$\supp \psi_f$  by
\begin{equation} 
\label{SuppDef}
    \supp \psi_f := \supp \omega_f,
    \quad
    \textrm { where } 
    \quad
    \omega_f (g) := \SP{\psi_f, \pi_\omega (g) \psi_f}.
\end{equation}
Here $f, g \in C^\infty_0 (M)[[\lambda]]$. Note that the 
$\mathbb C[[\lambda]]$-linear functional 
$\omega_f$ is positive and explicitly given by
\begin{equation}
    \omega_f (g) = \omega (\cc f * g * f).
\end{equation}
Clearly $\omega_f$ depends only on the equivalence class $\psi_f$. A  
first characterization of this support is given by the following
lemma: 
\begin{lemma}
\label{SuppPsifLem}
Let $\omega: C^\infty_0 (M)[[\lambda]] \to \mathbb C[[\lambda]]$ be a
positive $\mathbb C[[\lambda]]$-linear functional,
$f, g \in C^\infty_0 (M)[[\lambda]]$, and
$\alpha \in \mathbb C[[\lambda]]$.
\begin{enumerate}
\item If $\supp f \cap \supp \omega = \emptyset$ 
      then $f \in \mathcal J_\omega$ and thus $\psi_f = 0$.
\item $\supp \psi_f \subseteq \supp f \cap \supp \omega$. 
\item $\supp (\psi_f + \psi_g) \subseteq 
      \supp \psi_f \cup \supp \psi_g$ and 
      $\supp \alpha\psi_f \subseteq \supp \psi_f$.
\end{enumerate}
\end{lemma}
\begin{proof}
Let $\supp f \cap \supp \omega = \emptyset$ then also 
$\supp (\cc f * f) \cap \supp \omega = \emptyset$ whence 
$\omega (\cc f * f) = 0$ due to Lemma~\ref{BasicSuppLem} which proves
the first part. To avoid trivialities assume 
$\supp f \cap \supp \omega \ne \emptyset$. Considering
$g \in C^\infty_0 (M \setminus (\supp f \cap \supp \omega))$ we find 
$\supp (\cc f * g * f) \cap \supp \omega = \emptyset$ whence 
$\omega_f (g) = 0$ again due to Lemma~\ref{BasicSuppLem}. Thus the
second part follows. Finally we observe that for 
$h \in C^\infty_0 (M\setminus(\supp \psi_f \cup \supp \psi_g))$ 
we have $\omega (\cc f * h * g) \cc{\omega (\cc f * h * g)} \le
\omega (\cc f * f) \omega (\cc g * \cc h * h * g) = 0$ by the
Cauchy-Schwarz inequality since 
$\supp h \cap \supp \psi_g = \emptyset$. Hence we have  
$\omega (\cc f * h * g) = 0$ and similar we find 
$\omega (\cc g * h * f) = 0$. Now the third part follows easily since
for such $h$ one finds
$\omega_{f+g} (h) = \omega_f (h) + \omega_g(h) 
+ \omega (\cc f * h * g) + \omega (\cc g * h * f) = 0$.
\end{proof}

Note that in general $\supp \psi_f$ is strictly smaller than 
$\supp f \cap \supp \omega$. Examples can easily be found e.g.\ 
in the Schr{\"o}dinger representations in Section~\ref{ExSec}.

As a first consequence of this lemma we observe that for any open set
$O \subseteq M$ the space
\begin{equation}
\label{HODef}
    \mathfrak H_\omega (O) := \{ \psi_f \in \mathfrak H_\omega \;|\;
                                 \supp \psi_f \subseteq O \}
\end{equation}
is a sub-module of $\mathfrak H_\omega$ and clearly 
$\mathfrak H_\omega (O) \subseteq \mathfrak H_\omega (O')$ for 
$O \subseteq O'$. Thus the net structure of the observable algebra
$C^\infty_0 (M)[[\lambda]]$ induces a net structure for the GNS
pre-Hilbert space $\mathfrak H_\omega$. Note that it may happen that 
$\mathfrak H_\omega (O) = \{0\}$ for certain open $O \subseteq M$.
In order to characterize this net of pre-Hilbert spaces 
$\mathfrak H_\omega (O)$ we need the following lemma.
\begin{lemma}
\label{NiceRepLem}
Let $\omega : C^\infty_0 (M)[[\lambda]] \to \mathbb C[[\lambda]]$ be a
positive $\mathbb C[[\lambda]]$-linear functional and let 
$\psi_f, \psi_g \in \mathfrak H_\omega$.
\begin{enumerate}
\item If $\supp \psi_f \subseteq O$ with $O\subseteq M$ open then
      there exists a 
      $\tilde f \in \mathcal A_0 (O) 
      \subseteq C^\infty_0 (O)[[\lambda]]$ 
      such that $\psi_f = \psi_{\tilde f}$, i.e.\ 
      $f - \tilde f \in \mathcal J_\omega$.
\item $\supp \psi_f \cap \supp \psi_g = \emptyset$ implies 
      $\SP{\psi_f,\psi_g} = 0$.
\end{enumerate}
\end{lemma}
\begin{proof}
Choose an open neighborhood $U$ of $M\setminus O$ such that 
$U \cap \supp \psi_f = \emptyset$. Then let $\chi_0 + \chi_1 = 1$ be a
smooth partition of unity subordinate to the cover $U \cup O = M$,
i.e.\ $\supp \chi_0 \subseteq U$ and $\supp \chi_1 \subseteq O$. Since 
$\SP {\psi_{\chi_0*f}, \psi_{\chi_0*f}} 
= \omega_f (\cc \chi_0 * \chi_0) = 0$ due to 
$\supp \chi_0 \cap \supp \psi_f = \emptyset$ and 
Lemma~\ref{BasicSuppLem} one obtains 
$\psi_{\chi_0 * f} = 0$ whence $\psi_f = \psi_{\chi_1 * f}$. Setting
$\tilde f = \chi_1 * f$ the first part follows since clearly
$\supp \tilde f \subseteq \supp \chi_1 \subseteq O$.
Secondly, choose two
open sets $O, O'$ such that $O \cap O' = \emptyset$ and
$\supp \psi_f \subseteq O$ as well as $\supp \psi_g \subseteq O'$. Due
to the first part we may thus assume without restrictions that 
$f \in \mathcal A_0 (O)$ and 
$g \in \mathcal A_0 (O')$. Then 
$\cc f * g= 0$ whence $\SP {\psi_f, \psi_g} = \omega (\cc f * g) = 0$.
\end{proof}

This lemma suggests to consider those functions in the Gel'fand ideal 
$\mathcal J_\omega$ having their support in some open set $O$. We
define 
\begin{equation}
\label{JODef}
    \mathcal J_\omega (O) := 
    \mathcal A_0 (O) \cap \mathcal J_\omega
\end{equation}
for open $O \subseteq M$. 
Here we use $\mathcal A_0 (O)$ instead of $C^\infty_0 (O)[[\lambda]]$
and this difference will be crucial for the next proposition. By
definition $\mathcal J_\omega (O) \subseteq \mathcal A_0 (O)$
holds. Hence we can compare $\mathfrak H_\omega (O)$ with the quotient
$\mathcal A_0 (O) \big/ \mathcal J_\omega (O)$ which carries
a natural pre-Hilbert structure induced by the restriction of $\omega$
to $\mathcal A_0 (O)$. The next proposition states that they
are canonically isomorphic.
\begin{proposition}
\label{HomegaNetProp}
Let $\omega: C^\infty_0 (M) [[\lambda]] \to \mathbb C[[\lambda]]$ be a
positive $\mathbb C[[\lambda]]$-linear functional and $O,O' \subseteq M$
open.
\begin{enumerate}
\item $\mathfrak H_\omega (O)$ is canonically isometric to 
      $\mathcal A_0 (O) \big/ \mathcal J_\omega (O)$.
\item If $O \cap O' = \emptyset$ then
      $\mathfrak H_\omega (O) \; \bot \;\mathfrak H_\omega (O')$.
\item If $O \cap \supp \omega = \emptyset$ 
      then $\mathfrak H_\omega (O) = \{0\}$.
\end{enumerate}
\end{proposition}
\begin{proof}
Denote by $[f] \in \mathcal A_0 (O) \big/ \mathcal J_\omega (O)$ the
equivalence class of $f \in \mathcal A_0 (O)$. Then the Hermitian
product  
of $[f]$ and $[g]$ is given by $\SP{[f], [g]} = \omega (\cc f * g)$
and the canonical isomorphism to $\mathfrak H_\omega (O)$ is given by 
$[f] \mapsto \psi_f$. Since 
$\supp \psi_f \subseteq \supp f \subseteq O$ this is
clearly well-defined and isometric, 
hence injective. The surjectivity follows from 
Lemma~\ref{NiceRepLem} (i). The second part follows directly from 
Lemma~\ref{NiceRepLem} (ii). Finally let 
$O \cap \supp \omega = \emptyset$ and 
$\psi_f \in \mathfrak H_\omega (O)$. Then we can assume 
$f \in C^\infty_0 (O)[[\lambda]]$ due to 
Lemma~\ref{NiceRepLem} (i) whence 
$\SP{\psi_f, \psi_f} = \omega (\cc f * f) = 0$ due 
to Lemma~\ref{BasicSuppLem}. Hence $\psi_f = 0$, proving the 
third part. 
\end{proof}

As a first application of this proposition we have the following
corollary concerning convex sums of positive functionals:
\begin{corollary}
\label{ConvexCor}
Let $\omega_1, \omega_2 : C^\infty_0 (M)[[\lambda]] \to \mathbb C[[\lambda]]$ 
be two positive $\mathbb C[[\lambda]]$-linear functionals and let 
$\alpha_1, \alpha_2 \in \mathbb R[[\lambda]]$ be positive and let 
$\omega := \alpha_1 \omega_1 + \alpha_2 \omega_2$.
\begin{enumerate}
\item $\supp \omega = \supp \omega_1 \cup \supp \omega_2$.
\item If for $i=1,2$ one has $\supp \omega_i \subset O_i$ with 
      $O_1 \cap O_2 = \emptyset$ and $O_1, O_2$ open, then
      \begin{equation} 
      \label{HomegaSplit}
          \mathfrak H_\omega = \mathfrak H_\omega (O_1 \cup O_2) =
          \mathfrak H_\omega (O_1) \oplus \mathfrak H_\omega (O_2)
      \end{equation}
      and canonically $\mathfrak H_\omega (O_i) \cong \mathfrak H_{\omega_i}$.
\end{enumerate}
\end{corollary}
\begin{proof}
For the first part the inclusion `$\subseteq$' 
follows from Lemma~\ref{BasicSuppLem}
hence consider `$\supseteq$'. Assume $\omega_i (h) \ne 0$ for some $i$
and some $h \in C^\infty_0 (O)$. Choosing some 
$\chi \in C^\infty_0 (M)$ with $\chi = 1$ in an open neighborhood of
$\supp h$ one has $h = \chi * h$ and thus by the Cauchy-Schwarz
inequality one finds $\omega_i (\cc h * h) > 0$. Then 
$\omega (\cc h * h) \ge \alpha_i \omega_i (\cc h * h) > 0$ shows 
$\omega (\cc h * h) \ne 0$. Thus the first part follows since
$\supp (\cc h * h) \subseteq \supp h$. 
Secondly, consider $O_1 \cup O_2$ then 
$\mathfrak H_\omega = \mathfrak H_\omega (O_1 \cup O_2)$ 
since 
$\supp \psi_f \subseteq \supp \omega \subseteq O_1 \cup O_2$ for all 
$\psi_f \in \mathfrak H_\omega$ by Lemma~\ref{SuppPsifLem} (ii). 
Furthermore we may assume for $\psi_f$ that 
$f \in \mathcal A_0 (O_1 \cup O_2)$ due to Lemma~\ref{NiceRepLem} (i).
Clearly any such $f$ can be uniquely written as $f = f_1 + f_2$ where
$\supp f_i \subseteq O_i$ since $O_1 \cap O_2 = \emptyset$. 
But then $\psi_f = \psi_{f_1} + \psi_{f_2}$ 
with $\supp \psi_{f_i} \subseteq O_i$. This induces the above direct sum.
Finally, the canonical isomorphisms 
$\mathfrak H_\omega (O_i) \cong \mathfrak H_{\omega_i}$ 
are induced by the identity map which follows again from
Lemma~\ref{NiceRepLem}.
\end{proof}

Since the GNS pre-Hilbert space $\mathfrak H_\omega$ inherits the net 
structure $\mathfrak H_\omega (O)$ from the topology of $M$ for any
positive $\mathbb C[[\lambda]]$-linear functional $\omega$ we are
particularly interested in endomorphisms of $\mathfrak H_\omega$
respecting this locality structure. This motivates the following
definition of \emph{local} operators. 
\begin{definition}
\label{LocalOpDef}
Let $\omega: C^\infty_0 (M)[[\lambda]] \to \mathbb C[[\lambda]]$ be a
positive $\mathbb C[[\lambda]]$-linear functional and 
$A: \mathfrak H_\omega \to \mathfrak H_\omega$ a 
$\mathbb C[[\lambda]]$-linear map. Then $A$ is called local if for all
$\psi_f \in \mathfrak H_\omega$
\begin{equation}
\label{LocalOperatorDef}
    \supp A\psi_f \subseteq \supp \psi_f.
\end{equation}
The set of local $\mathbb C[[\lambda]]$-linear endomorphisms of 
$\mathfrak H_\omega$ is denoted by $\Loc (\mathfrak H_\omega)$.
\end{definition}
Similarly one defines the local $\mathbb C[[\lambda]]$-linear maps 
$\Loc (\mathfrak H_\omega, \mathfrak H_{\omega'})$ from one GNS
pre-Hilbert space $\mathfrak H_\omega$ into another 
$\mathfrak H_{\omega'}$, where $\omega$ and $\omega'$ are two positive
$\mathbb C[[\lambda]]$-linear functionals. The following proposition
is an obvious consequence of the preceding lemmas. 
\begin{proposition}
Let $\omega: C^\infty_0 (M)[[\lambda]] \to \mathbb C[[\lambda]]$ be a
positive $\mathbb C[[\lambda]]$-linear functional then 
$\Loc (\mathfrak H_\omega)$ is a subalgebra of all endomorphisms of
$\mathfrak H_\omega$. 
\end{proposition}
The abstract notion of local operators in GNS pre-Hilbert spaces
turns out to be fundamental for the whole following discussion. In
fact the GNS representation always is local.
\begin{theorem}
Let $(M, *)$ be a Poisson manifold with local star product and let 
$\omega: C^\infty_0 (M)[[\lambda]] \to \mathbb C[[\lambda]]$
be a positive $\mathbb C[[\lambda]]$-linear functional with GNS
pre-Hilbert space $\mathfrak H_\omega$. Then for any 
$f \in C^\infty (M)[[\lambda]]$ the operator 
$\pi_\omega (f): \mathfrak H_\omega \to \mathfrak H_\omega$ is
local. Moreover, 
\begin{equation}
\label{LocalGNSOp}
    \supp \pi_\omega (f)\psi_g \subseteq \supp f \cap \supp \psi_g
\end{equation}
for all $f \in C^\infty (M)[[\lambda]]$ and 
$\psi_g \in \mathfrak H_\omega$. 
\end{theorem}
\begin{proof}
We have to determine the support of $\psi_{f*g}$. Without restriction 
assume $\supp \psi_g \ne M$ and consider 
$h \in C^\infty_0 (M \setminus \supp \psi_g)$. Then  
$\SP{\psi_{f*g}, \pi_\omega (h) \psi_{f*g}} 
= \omega (\cc g * \cc f * h * f *g) = \omega_g (\cc f * h * f) = 0$ 
since $\supp (\cc f * h * f) \subseteq \supp h$. Thus 
$\supp \psi_{f*g} \subseteq \psi_g$ follows. Moreover, 
$\supp \pi_\omega (f)\psi_g \subseteq 
\supp \omega \cap \supp g \cap \supp f$ due to
Lemma~\ref{SuppPsifLem} (ii) and the locality of $*$. Thus
(\ref{LocalGNSOp}) follows.
\end{proof}

We shall now turn to the question how these local operators can be
used to give a reasonable distinction between the various types of GNS 
representations. In particular we are searching for a suitable way to
distinguish `thermal' representations from `pure' ones. In the usual
$C^*$-algebra theory as used e.g.\ in quantum mechanics and quantum
field theory the concept of `thermal' and `pure' states is highly
developed both from the physical and mathematical point of view and
can be summarized as follows. Pure states are the extremal points in
the convex set of all states and a state is pure if and only if its
GNS representation is irreducible which is the case if and only if the
commutant in the GNS representation is trivial, see
e.g.\ \cite[Thm.~2.3.19]{BR87}. On the other hand the so-called
KMS states, which
are understood to describe `thermal' behavior, are not pure states. In
deformation quantization it will turn out that the commutant of the
GNS representation within the \emph{local operators} is the
appropriate choice for the definition of a commutant which
distinguishes `thermal' from `pure' representations.

Before we investigate local commutants let us firstly consider the
adjoints of local operators.
\begin{lemma}
Let $\omega: C^\infty_0 (M)[[\lambda]] \to \mathbb C[[\lambda]]$ be a
positive $\mathbb C[[\lambda]]$-linear functional with GNS pre-Hilbert
space $\mathfrak H_\omega$. If $A \in \Loc (\mathfrak H_\omega)$ has
an adjoint operator $A^*$ then $A^*$ is local again.
\end{lemma}
\begin{proof}
Assume $A^*$ exists and let $\psi_f \in \mathfrak H_\omega$ be
arbitrary where we may assume $\supp \psi_f \ne M$ to avoid
trivialities.  
Choose $h \in C^\infty_0 (M)[[\lambda]]$ such that 
$A^* \psi_f = \psi_h$. Then it is sufficient to show 
$\omega_h |_{C^\infty_0 (M \setminus \supp \psi_f)} = 0$ since in
this case $\supp A^*\psi_f = \supp \psi_h \subseteq \supp \psi_f$. 
Hence let $g \in C^\infty_0 (M \setminus \supp \psi_f)$
then we compute $\omega_h (g)$ using $(A^*)^* = A$ and 
$\pi_\omega (g)^* = \pi_\omega (\cc g)$ and obtain
$\omega_h (g) = \SP{A\pi_\omega (\cc g) A^*\psi_f, \psi_f}$.
Now $\pi_\omega (\cc g)$ is local and due to (\ref{LocalGNSOp}) we
even have $\supp \pi_\omega (\cc g) A^* \psi_f \subseteq \supp g$.
With  $\supp g \cap \supp \psi_f = \emptyset$ we conclude that
$\omega_h (g) = 0$ using Lemma~\ref{BasicSuppLem}. Thus $A^*$ is
local. 
\end{proof}

It will turn out that in various examples the local operators
automatically have an adjoint though this is not evident from the
abstract point of view. We now define
\begin{equation}
\label{AdLocalOpDef}
    \BLoc (\mathfrak H_\omega) := \Loc (\mathfrak H_\omega) \cap 
                                 \Bounded (\mathfrak H_\omega),
\end{equation}
and conclude from the last lemma that $\BLoc (\mathfrak H_\omega)$ is
a $^*$-algebra.
\begin{proposition}
Let $\omega: C^\infty_0 (M)[[\lambda]] \to \mathbb C[[\lambda]]$ be a
positive $\mathbb C[[\lambda]]$-linear functional with GNS pre-Hilbert
space $\mathfrak H_\omega$. Then $\BLoc (\mathfrak H_\omega)$ is a
$^*$-algebra over $\mathbb C[[\lambda]]$.
\end{proposition}

Using the algebras $\Loc (\mathfrak H_\omega)$ and
$\BLoc (\mathfrak H_\omega)$, respectively, we can define the
commutant of an arbitrary subset of $\Loc (\mathfrak H_\omega)$ and
$\BLoc (\mathfrak H_\omega)$, respectively, in the following way:
\begin{definition}
Let $\omega: C^\infty_0 (M)[[\lambda]] \to \mathbb C[[\lambda]]$ be a
positive $\mathbb C[[\lambda]]$-linear functional with GNS pre-Hilbert
space  $\mathfrak H_\omega$ and let 
$\mathcal B \subseteq \Loc (\mathfrak H_\omega)$ 
or $\mathcal B \subseteq \BLoc (\mathfrak H_\omega)$ then
\begin{equation}
    \begin{array} {rcl}
        \mathcal B'_{\Loc} 
        & := & \{ A \in \Loc (\mathfrak H_\omega) \; | \;
               \forall B \in \mathcal B: \; AB = BA \}, \\
        \mathcal B'_{\BLoc}
        & := & \{ A \in \BLoc (\mathfrak H_\omega) \; | \; 
               \forall B \in \mathcal B: \; AB = BA \}
    \end{array}
\end{equation}
are called the local and $^*$-local commutant of $\mathcal B$ in 
$\Loc (\mathfrak H_\omega)$ and $\BLoc (\mathfrak H_\omega)$,
respectively. 
\end{definition}

In principle one could also define a $\Bounded$-commutant in 
$\Bounded (\mathfrak H_\omega)$ but this will not be as useful as the
above versions of commutants. If the context is clear we will
sometimes omit the subscript $\Loc$ resp. $\BLoc$. Now let 
$\mathcal B \subseteq \Loc (\mathfrak H_\omega)$ or
$\BLoc (\mathfrak H_\omega)$, respectively. Then $\mathcal B'$ is a
subalgebra of $\Loc (\mathfrak H_\omega)$ or 
$\BLoc (\mathfrak H_\omega)$, respectively, and if 
$\mathcal B = \mathcal B^* \subseteq \BLoc (\mathfrak H_\omega)$ then 
$\mathcal B'_{\BLoc}$ is even a $^*$-subalgebra. Note that 
$\mathcal B$ needs not to be an algebra at all. Furthermore one has 
$\mathcal B \subseteq \mathcal B''$ and if 
$\mathcal C \subseteq \mathcal B$ then 
$\mathcal B' \subseteq \mathcal C'$. Hence 
$\mathcal B''' = \mathcal B'$ and thus the commutant of a set is
always closed under taking the double commutant. Note finally that for
a subalgebra $\mathcal B$ the algebra 
$\mathcal Z = \mathcal B \cap \mathcal B'$ is the center of 
$\mathcal B$ and $\mathcal B'$. As for complex Hilbert spaces 
we shall call a
$^*$-subalgebra of $\BLoc (\mathfrak H_\omega)$ a 
\emph{von Neumann algebra} if  
$\mathcal B = \mathcal B''_{\BLoc}$ and a \emph{factor} if in addition
the center is trivial. Here we are forced to take $\BLoc$-commutants
since we are interested in $^*$-algebras.

Now we can use this notion of local commutants to characterize `pure' 
and `thermal' GNS representations. A first result is obtained in the  
following proposition which can be viewed as an analogue to the
well-known situation of mixed states for complex $C^*$-algebras, see
e.g.\ \cite[Thm.~2.3.19]{BR87}.
\begin{proposition}
Let 
$\omega_1, \omega_2: C^\infty_0 (M)[[\lambda]] \to \mathbb C[[\lambda]]$
be two non-zero positive $\mathbb C[[\lambda]]$-linear functionals with 
$\supp \omega_1 \cap \supp \omega_2 = \emptyset$. Then the local
(and $^*$-local) commutant of $\pi_\omega (\mathcal A)$ is non-trivial,
where $\omega = \alpha_1 \omega_1 + \alpha_2 \omega_2$, provided 
$\alpha_1, \alpha_2 > 0$.
\end{proposition}
\begin{proof}
Since $\supp \omega_1 \cap \supp \omega_2 = \emptyset$ we can find two 
open sets $O_1, O_2$ such that $\supp \omega_i \subseteq O_i$ and 
$O_1 \cap O_2 = \emptyset$. Now we apply Cor.~\ref{ConvexCor} to 
obtain 
$\mathfrak H_\omega = 
\mathfrak H_\omega (O_1) \oplus \mathfrak H_\omega (O_2)$.
Then clearly the projectors on $\mathfrak H_\omega (O_i)$, $i=1,2$,
commute with $\pi_\omega (f)$ for all $f \in \mathcal A$ and are
clearly $^*$-local operators proving the proposition.
\end{proof}

Note that with the hypothesis of the above proposition 
any function $\chi \in C^\infty (M)[[\lambda]]$ with 
$\chi|_{O_i} = c_i$, where $c_i \in \mathbb C[[\lambda]]$ are
constants, is in the center of $\pi_\omega (C^\infty (M)[[\lambda]])$ 
but acts non-trivial on $\mathfrak H_\omega$ if $c_1 \ne c_2$. 
\begin{corollary}
\label{SuperSelectionCor}
Let $\omega: C^\infty_0 (M)[[\lambda]] \to \mathbb C[[\lambda]]$ be a
positive $\mathbb C[[\lambda]]$-linear functional such that
$\supp \omega$ has at least two connected components. Then 
$\pi_\omega (C^\infty (M)[[\lambda]])$ has a non-trivial local
(and $^*$-local) commutant. 
\end{corollary}

%
%

\section{Faithful positive linear functionals}
\label{FaithSec}

Let us now consider faithful positive $\mathbb C[[\lambda]]$-linear 
functionals and their GNS representations. It turns out that they can 
completely be characterized by their support.

First recall that a positive $\mathbb C[[\lambda]]$-linear functional
$\omega$ is called \emph{faithful} if $\mathcal J_\omega = \{0\}$. 
Hence the GNS pre-Hilbert space $\mathfrak H_\omega$ is
canonically isomorphic to $C^\infty_0 (M)[[\lambda]]$ via 
$\mathfrak H_\omega \ni \psi_f 
\mapsto f \in C^\infty_0 (M)[[\lambda]]$ 
as $\mathbb C[[\lambda]]$-modules. It will sometimes
be useful not to identify $\mathfrak H_\omega$ and 
$C^\infty_0 (M)[[\lambda]]$ but use this isomorphism since 
$\mathfrak H_\omega$ has the Hermitian
product as additional structure. Under the above isomorphism the
corresponding GNS representation is simply given by left
multiplication in $C^\infty_0 (M)[[\lambda]]$. Thus we use also
the notion $\mathsf L_f$ instead of $\pi_\omega (f)$ for this particular
representation where $f \in C^\infty (M)[[\lambda]]$. On the other
hand a representation $\pi$ is called faithful if it is injective.

The following technical lemma concerning local left inverses is proved
by the usual recursion techniques.
\begin{lemma}
Let $(M,*)$ be a Poisson manifold with local star product,
$f = \sum_{r=0}^\infty \lambda^r f_r$ with $f_r \in C^\infty_0 (O)$,
and $\emptyset \ne U \subseteq O$ an open subset such that $f_0 (x) \ne 0$ 
for all $x \in U$. Then there exists another non-empty open subset 
$U' \subseteq U$ and a function $f^{-1} \in C^\infty_0 (U)[[\lambda]]$
such that $f^{-1} * f |_{U'} = 1$. 
\end{lemma}
Similarly there exist local right inverses such that 
$f * \tilde f^{-1}|_{U''} = 1$ and on $U' \cap U''$ left and right 
inverses coincide. Using such a local left inverse the following
proposition can be shown easily.
\begin{proposition}
\label{FaithfulSuppProp}
Let $\omega: C^\infty_0 (M)[[\lambda]] \to \mathbb C[[\lambda]]$ be a
positive $\mathbb C[[\lambda]]$-linear functional. Then $\omega$ is
faithful if and only if $\supp \omega = M$. 
\end{proposition}
\begin{proof}
Let $\omega$ be faithful and let 
$0 \ne f \in C^\infty_0 (O)$ for some non-empty open set $O$. Then 
$\supp (\cc f * f) \subseteq O$ and thus 
$\omega (\cc f * f) > 0$. This implies $\supp \omega = M$. Now assume 
$\supp \omega = M$, and assume we have found a function 
$0 \ne f \in C^\infty_0 (O)[[\lambda]]$ such that 
$\omega (\cc f * f) = 0$. Without  restriction we can assume that
already the lowest order of $f$ is non-zero. Hence there is a
non-empty open subset $U$ and a local left inverse  
$f^{-1} \in C^\infty_0 (O)[[\lambda]]$ with 
$f^{-1} * f |_{U} = 1$. Now let $h \in C^\infty_0 (U)[[\lambda]]$ 
be arbitrary then clearly $h = g * f$ with some 
$g \in C^\infty_0 (U)[[\lambda]]$, namely $g = h * f^{-1}$. Then 
$\omega(h) \cc{\omega(h)} \le \omega(\cc g *g) \omega(\cc f * f) = 0$
shows $\omega (h) = 0$ for all $h \in C^\infty_0 (U)[[\lambda]]$. Thus 
$\omega|_{C^\infty_0 (U)[[\lambda]]} = 0$ in contradiction to 
$\supp \omega = M$.
\end{proof}

Let us now investigate the support of $\psi_f \in \mathfrak H_\omega$
for a faithful positive functional $\omega$. Since in this case 
canonically $\mathfrak H_\omega \cong C^\infty_0 (M)[[\lambda]]$ as 
$\mathbb C[[\lambda]]$-modules we expect that the support of $\psi_f$
coincides with the support of $f$. This is indeed the case.
\begin{lemma}
Let $\omega: C^\infty_0 (M)[[\lambda]] \to \mathbb C[[\lambda]]$ be a
faithful positive $\mathbb C[[\lambda]]$-linear functional. Then for
all $f \in C^\infty_0 (M)[[\lambda]]$ we have 
$\supp \psi_f = \supp f$. 
\end{lemma}
\begin{proof}
The inclusion $\subseteq$ is in general true due to 
Lemma~\ref{SuppPsifLem} (ii). Hence let 
$f \in C^\infty_0(M)[[\lambda]]$ and assume $\supp \psi_f \ne M$ to
avoid trivialities. Then for all 
$g \in C^\infty_0 (M \setminus \supp \psi_f)$ we have 
$0 = \omega_f (\cc g * g) = \omega (\cc{(g*f)} * (g*f))$ whence 
$g*f = 0$ since $\omega$ is faithful. But this implies 
$\supp f \cap (M \setminus\supp\psi_f) = \emptyset$ and thus the claim 
follows.
\end{proof}

Since a faithful positive functional $\omega$ has Gel'fand ideal 
$\mathcal J_\omega = \{0\}$ the \emph{right multiplication} 
$\mathsf R_f$ by $f \in C^\infty (M)[[\lambda]]$ 
\begin{equation}
\label{RightMult}
    \mathsf R_f \psi_g := \psi_{g*f}, \qquad \psi_g \in \mathfrak H_\omega
\end{equation}
is well-defined for all $f \in C^\infty (M)[[\lambda]]$ and clearly
again a local operator. 
Note that in general $\mathsf R_f$ is well-defined if $f$ is contained
in the \emph{Lie idealizer} of the Gel'fand ideal, i.e.\ the largest
subalgebra of $C^\infty (M)[[\lambda]]$ containing 
$\mathcal J_\omega$ as a two-sided ideal. Using this right
multiplication and the preceding results we obtain immediately the
following corollaries:
\begin{corollary}
Let $\omega: C^\infty_0 (M)[[\lambda]] \to \mathbb C[[\lambda]]$ be a
faithful positive $\mathbb C[[\lambda]]$-linear functional. Then an
operator $A: \mathfrak H_\omega \to \mathfrak H_\omega$ is local if
and only if the corresponding operator 
$A: C^\infty_0 (M)[[\lambda]] \to C^\infty_0 (M)[[\lambda]]$ is local
in the usual sense. 
\end{corollary}
\begin{corollary}
\label{FaithfulRepCommCor}
Let $\omega: C^\infty_0 (M)[[\lambda]] \to \mathbb C[[\lambda]]$ be a 
faithful positive $\mathbb C[[\lambda]]$-linear functional. Then the
local commutant of $\pi_\omega (C^\infty (M)[[\lambda]])$ contains all
right multiplications $\mathsf R_f$ with 
$f \in C^\infty(M)[[\lambda]]$. Thus  
$\pi_\omega (C^\infty (M)[[\lambda]])_{\Loc}'$ is non-trivial 
(if $\dim M > 0$). 
\end{corollary}

The question whether the $^*$-local commutant is non-trivial seems to
be more complicated since $\supp \omega = M$ still allows rather
`wild' functionals whence existence of an adjoint of $\mathsf R_f$ is
not obvious. Take e.g.\ in zeroth order a faithful positive linear
functional, as e.g.\ integration over some positive density, then one
can add in higher orders of $\lambda$ \emph{any} real linear
functionals and still has a positive linear functional. 
Nevertheless things become simpler if we consider 
KMS functionals in Section~\ref{ExSec}.

Let us finally consider the GNS representation of such faithful
functionals: 
\begin{proposition}
Let $\omega: C^\infty_0 (M)[[\lambda]] \to \mathbb C[[\lambda]]$ be a
positive $\mathbb C[[\lambda]]$-linear functional with GNS representation
$\pi_\omega$. Then $\pi_\omega$ is faithful if and only if $\omega$ is
faithful. 
\end{proposition}
\begin{proof}
Note that $\pi_\omega$ is always understood to be extended to
$C^\infty (M)[[\lambda]]$. Assume first that $\omega$ is faithful then
$C^\infty (M)[[\lambda]]$ is represented on 
$\mathfrak H_\omega \cong C^\infty_0 (M)[[\lambda]]$ by left
multiplications which is clearly faithful. On the other hand assume
that $\omega$ is not faithful. Then $O = M \setminus \supp \omega$ is
a non-empty open subset due to Prop.~\ref{FaithfulSuppProp}. Due to
(\ref{LocalGNSOp}) we have $\pi_\omega (f) = 0$ for all 
$f \in C^\infty_0 (O)$. Thus $\pi_\omega$ cannot be faithful.
\end{proof}

We observe that if a GNS representation is faithful then it is 
equivalent to the left multiplication of elements of 
$C^\infty (M)[[\lambda]]$ on $C^\infty_0 (M)[[\lambda]]$. Hence a
faithful GNS representation $\pi_\omega$ does not really depend on
$\omega$ but is uniquely given. Nevertheless the Hermitian product of 
$\mathfrak H_\omega$ still depends crucially on $\omega$. Moreover, a
faithful GNS representation has always a non-trivial local
commutant. Thus in this point the situation is quite different from
the usual $C^*$-algebra theory where a faithful GNS representation can
of course be irreducible, as e.g.\ the standard  representation of the
bounded operators on a Hilbert space.

%
%

\section{Basic examples}
\label{ExSec}

\subsection*{Traces and KMS functionals}

As first basic example we consider positive traces and KMS
functionals. A trace $\tr$ of the algebra $C^\infty (M)[[\lambda]]$ is
a $\mathbb C[[\lambda]]$-linear functional defined on 
$C^\infty_0 (M)[[\lambda]]$ such that $\tr (f * g) = \tr (g * f)$. 
In case where $M$ is \emph{symplectic and connected} there exists up
to normalization a unique such trace \cite{NT95a,NT95b}. Moreover,
this trace functional is of the form 
\begin{equation}
\label{TraceFunct}
    \tr (f) = c \int_M \left(f + \sum_{r=1}^\infty 
              \lambda^r T_r (f) \right) \Omega
\end{equation}
where $\Omega = \omega \wedge \cdots \wedge \omega$ is the Liouville
form, $c \in \mathbb C[[\lambda]]$ a normalization factor, and the
$T_r$ are differential operators. If the star product satisfies 
$\cc{f  * g} = \cc g * \cc f$, as we assume,  then it can easily be
shown that by an appropriate choice of the normalization factor $c$
the trace becomes a real functional, i.e.\ 
$\tr (\cc f) = \cc{\tr (f)}$. Since in lowest order of $\lambda$ the
trace then consists in integration over $M$ it follows from
\cite[Lem.~2]{BW98a} that $ \tr$ becomes a positive functional
\cite[Lem.~4.3]{BRW98a}.

Given such a positive trace $\tr$ it is easily seen that 
$\supp \tr = M$ and thus 
the GNS representation of $\tr$ is faithful according to the last
section. Hence 
$\mathfrak H_{\tr} \cong C^\infty_0 (M)[[\lambda]]$ and $\pi_\tr$ is
equivalent to  
the left multiplication $\mathsf L$. Moreover the Hermitian product is
given by $\SP{\psi_f, \psi_g} = \tr(\cc f * g)$. It follows that not
only the local commutant is non-trivial as stated in 
Corollary~\ref{FaithfulRepCommCor} but even the $^*$-local commutant is
non-trivial: consider an arbitrary element 
$f \in C^\infty (M)[[\lambda]]$ and the corresponding right
multiplication operator $\mathsf R_f$ on $\mathfrak H_{\tr}$. Clearly 
$\mathsf R_f \in (\pi_{\tr} (C^\infty (M)[[\lambda]]))'_{\Loc}$ but
now we can even prove the existence of $\mathsf R_f^*$ whence 
$\mathsf R_f \in (\pi_{\tr} (C^\infty (M)[[\lambda]]))'_{\BLoc}$ for every 
$f \in C^\infty (M)[[\lambda]]$. Namely $\mathsf R_f^*$ is given by
\begin{equation}
    \mathsf R_f^* = \mathsf R_{\cc f},
\end{equation}
as an easy computation shows. 
Thus $(\pi_\tr (C^\infty (M)[[\lambda]]))'_{\BLoc}$ is non-trivial
too.

Another important example is given by the so-called KMS functionals.  
In \cite{BFLS84,BL85,BRW98a,BRW98b} the notion of 
KMS states known from $C^*$-algebra theory, see 
e.g.\ \cite{Con94,BR81,BR87,Haa93}, was transfered to the framework of 
deformation quantization. We shall only use the final result on the 
existence and uniqueness of these functionals and their particular form 
as found in \cite{BRW98a,BRW98b}. To this end we first need a notion
of star exponential (see e.g.\ \cite{BFFLS78}), i.e.\ the analogue of the
exponential series build out of star product powers. For our purpose
it is sufficient to use the following definition avoiding questions on 
convergence. The \emph{star exponential}
$\Exp(\beta H) \in C^\infty (M)[[\lambda]]$ of 
$H \in C^\infty (M)[[\lambda]]$ with $\beta \in \mathbb R$ is defined
to be the unique solution in $C^\infty (M)[[\lambda]]$ of the
differential equation  
\begin{equation}
\label {StarExpDef}
    \frac{d}{d\beta} \Exp (\beta H) = H * \Exp (\beta H)
\end{equation}
with initial condition $\Exp (0) = 1$. Of course one has to show the 
existence and uniqueness of such a solution but this has been done
e.g.\ in \cite[Lem.~2.2]{BRW98a}. Moreover, the usual properties hold,
i.e.\ $\Exp(\beta H)$ commutes with $H$ and satisfies 
$\Exp((\beta+\beta')H) = \Exp (\beta H) * \Exp (\beta' H)$ as well as
$\cc{\Exp (\beta H)} = \Exp (\beta \cc H)$ for all 
$\beta,\beta' \in \mathbb R$ and $H \in C^\infty (M)[[\lambda]]$ (due
to $\cc{f*g} = \cc g * \cc f$). Note finally, that $\Exp (\beta H)$
can even be extended to arbitrary $\beta \in \mathbb C[[\lambda]]$
such that the above relations hold.
Then in \cite{BRW98a} it was shown that
for a given real `Hamiltonian' $H \in C^\infty (M)[[\lambda]]$ and a
given `inverse temperature' $\beta \in \mathbb R$ there exists an up to
normalization unique KMS functional 
$\omega_\KMS: C^\infty_0 (M)[[\lambda]] \to \mathbb C[[\lambda]]$
given by  
\begin{equation}
\label{KMSFunctional}
    \omega_\KMS (f) = \tr (\Exp (-\beta H)*f),
\end{equation}
where $\tr$ is the trace of $C^\infty (M)[[\lambda]]$. We observe that
for a \emph{positive} trace the KMS functional (\ref{KMSFunctional}) is
positive too, since $\cc H = H$. Moreover, 
$\supp \omega_\KMS = M$ since $\Exp (-\beta H)$ 
is invertible. Hence again 
$\mathfrak H_\KMS \cong C^\infty_0 (M)[[\lambda]]$  
as $\mathbb C[[\lambda]]$-module and the corresponding GNS
representation is again given by left multiplication. Now the
Hermitian product is given by
\begin{equation}
\label{KMSHermiteanProd}
    \SP{\psi_f, \psi_g}_\KMS
    = \tr (\Exp (-\beta H) * \cc f * g).
\end{equation}
Again we see that the right multiplication $\mathsf R_f$ is contained in 
$(\pi_\KMS (C^\infty (M)[[\lambda]]))'_{\Loc}$. A straightforward
computation shows that $\mathsf R_f^*$ exists and is given by
$\mathsf R_f^* = \mathsf R_{\Exp (-\beta H) * \cc f * \Exp (\beta H)}$.
Thus $\mathsf R_f \in (\pi_\KMS (C^\infty (M)[[\lambda]]))'_{\BLoc}$
for all $f \in C^\infty (M)[[\lambda]]$ and hence even the $^*$-local
commutant of the GNS representation is non-trivial for KMS
functionals. This result is expected for physical reasons since KMS
functionals are believed to describe \emph{thermal behavior} and physical
situations in thermal equilibrium whence they should be `mixed'. 
Note that $\beta = 0$ (i.e.\ infinite temperature) brings us back 
to the case of the
positive trace. We summarize these results in a proposition: 
\begin{proposition}
Let $(M, *)$ be a connected symplectic manifold with local star
product and positive trace $\tr$. Let 
$H \in C^\infty (M)[[\lambda]]$ be a
real Hamiltonian and $\beta \in \mathbb R$. Denote by $\omega_\KMS$
the corresponding positive KMS functional with GNS pre-Hilbert space
$\mathfrak H_\KMS$. Then one has: 
\begin{enumerate}
\item $\supp \omega_\KMS = M$ 
      whence $\mathfrak H_\KMS \cong C^\infty_0 (M)[[\lambda]]$.
\item $\mathsf R_f \in \BLoc (\mathfrak H_\KMS)$ for all 
      $f \in C^\infty (M)[[\lambda]]$  and
      \begin{equation}
      \label{KMSAdjoint}
          \mathsf R_f^* =
          \mathsf R_{\Exp(-\beta H) * \cc f * \Exp (\beta H)}.
      \end{equation}
\item $\mathsf R_f \in (\pi_\KMS (C^\infty (M)[[\lambda]]))'_{\BLoc}$
      for all $f \in C^\infty (M)[[\lambda]]$ and thus the $^*$-local
      commutant of the GNS representation is non-trivial (if $\dim M >0$).
\end{enumerate}
\end{proposition}

Let us finally investigate the relation between two KMS functionals
and their GNS representations. Let $\omega_\KMS$ and 
$\omega'_\KMS$ be the positive KMS functionals for 
$(H, \beta)$ and $(H',\beta')$, respectively, 
normalized in the same way (\ref{KMSFunctional}). It turns out that
the GNS representations are locally and unitarily equivalent
by an explicitly given unitary map. Remember that this fact is quite
different from the situation in quantum field theory where the GNS
representations of KMS states for different temperatures are known to
be unitarily inequivalent under quite general pre-conditions, see
e.g.\ \cite{Tak70}. The main point is that the usual representations
are type III representations \cite{Haa93} 
and thus our result suggests that deformation quantization somehow
corresponds not to a type~III representation. The trace being defined
on a twosided ideal $C^\infty_0 (M)[[\lambda]]$ of the whole algebra
$C^\infty (M)[[\lambda]]$ reminds much more on a type I
representation, if heuristically deformation quantization is interpreted
as an asymptotic expansion for $\hbar \to 0$ of some convergent
situation. This is not surprising since star products on
\emph{finite-dimensional} symplectic manifolds correspond physically to
a finite number of degrees of freedom. If $M$ is even compact, this
even looks like a matrix algebra and hence a type I$_n$
with $n \in \mathbb N$. Thus a compact symplectic manifold somehow
corresponds to a finite dimensional Hilbert space, a result which is
also obtained in other quantizations schemes as e.g.\ geometric
quantization \cite{Woo92}. Let us now state the result precisely: 
\begin{proposition}
\label{AllKMSEquivProp}
Let $(M, *)$ be a connected symplectic manifold with local star
product and positive trace $\tr$. Then for any two real Hamiltonians 
$H, H' \in C^\infty (M)[[\lambda]]$ and any two inverse temperatures 
$\beta, \beta' \in \mathbb R$ the GNS representations 
$\pi_\KMS$ and $\pi'_\KMS$ of the corresponding KMS
functionals $\omega_\KMS$ and $\omega'_\KMS$
are locally unitarily equivalent via the unitary map 
$U : \mathfrak H_\KMS \to \mathfrak H'_\KMS$ given by 
\begin{equation}
\label{UniKMSEquiv}
    \mathfrak H_\KMS \ni \psi_f \; \mapsto \; U \psi_f = 
    \psi_{f * \Exp (-\frac{\beta}{2} H) * \Exp (\frac{\beta'}{2}H')}
    \in \mathfrak H'_\KMS.
\end{equation}
\end{proposition}
\begin{proof}
Notice that $U$ is indeed well-defined and clearly 
$U \in \Loc(\mathfrak H_\KMS, \mathfrak H'_\KMS)$. Then the unitary
equivalence is a simple computation.
\end{proof}

\subsection*{Bargmann-Fock representation on K{\"a}hler manifolds}

Another fundamental example is given by the formal Bargmann-Fock 
representation on K{\"a}hler manifolds. Before we discuss the general 
situation let us briefly remember the well-known situation for 
$M = \mathbb C^n$. Viewing $\mathbb C^n$ as K{\"a}hler manifold
with global holomorphic coordinates $z^1, \ldots, z^n$ and endowed
with the usual symplectic (K{\"a}hler) form  
$\omega = \frac{\im}{2} \sum_{k=1}^n dz^k \wedge d\cc z^k$ we
consider the \emph{Wick star product}
\begin{equation}
\label{WickStarProd}
    f \wick g = \sum_{r=0}^\infty \frac{(2\lambda)^r}{r!} 
            \sum_{i_1,\ldots,i_r} 
            \frac{\partial^r f}
            {\partial z^{i_1} \cdots \partial z^{i_r}} 
            \frac{\partial^r g}
            {\partial \cc z^{i_1} \cdots \partial \cc z^{i_r}},
\end{equation}
where $f, g \in C^\infty (\mathbb C^n)[[\lambda]]$, see e.g.\ 
\cite{BW97a}. It
turns out that the evaluation functional $\delta_p$ at any point 
$p \in \mathbb C^n$ is positive with respect to $\wick$, and clearly
the support of $\delta_p$ is given by $\{p\}$. Considering for
simplicity the point $p=0$ one finds that
the Gel'fand ideal of $\delta_0$ is given by \cite[Lem.~7]{BW98a}
\begin{equation}
\label{WickGelfand}
    \mathcal J_0 = \left\{ f \in C^\infty_0 (\mathbb C^n)[[\lambda]]
                   \; \left| \; 
                   \forall I: \;
                   \frac{\partial^{|I|} f}
                   {\partial \cc z^I} (0) = 0 \right.\right\},
\end{equation}
where $I = (\cc i_1, \ldots, \cc i_r)$, $r \ge 0$, ranges over all
multi-indices. Finally, one obtains that the GNS pre-Hilbert space
$\mathfrak H_0$ can be described by 
$(\mathbb C[[\cc y^1, \ldots, \cc y^n]])[[\lambda]]$ where 
\begin{equation}
\label{WickIso}
    \mathfrak H_0 \ni \psi_f \mapsto 
    \sum_{r=0}^\infty \sum_{i_1, \ldots, i_r} \frac{1}{r!}
    \frac{\partial^r f}
    {\partial \cc z^{i_1} \cdots \partial \cc z^{i_r}} (0)\,
    \cc y^{i_1} \cdots \cc y^{i_r}
\end{equation}
is the isomorphism. Since $\supp \delta_0 = \{0\}$ is only a single
point we see that $\supp \psi_f = \{0\}$ for any $\psi_f \ne 0$. Thus
\emph{any} $\mathbb C[[\lambda]]$-linear endomorphism of
$\mathfrak H_0$ is necessarily local. In order to compute the local
commutant of $\pi_0 (C^\infty (\mathbb C^n)[[\lambda]])$ we recall 
from \cite[Lem.~8]{BW98a} that the GNS representation is given by the
formal analogue of the \emph{Bargmann-Fock representation}
\begin{equation}
\label{BargmannFock}
    \pi_0 (f) = \sum_{r,s=0}^\infty \frac{(2\lambda)^r}{r!s!}
                \sum_{{i_1, \ldots, i_r \atop j_1, \ldots, j_s}}
                \frac{\partial^{r+s} f}
                {\partial z^{i_1} \cdots \partial z^{i_r}
                 \partial \cc z^{j_1} \cdots \partial \cc z^{j_s}}
                (0)\,
                \cc y^{j_1} \cdots \cc y^{j_s}
                \frac{\partial^r}
                {\partial \cc y^{i_1} \cdots \cc y^{i_r}},
\end{equation}
where we used the isomorphism (\ref {WickIso}). Then we obtain the
following result:
\begin{proposition}
\label{WickCnProp}
Let 
$\delta_0: C^\infty_0 (\mathbb C^n)[[\lambda]] \to \mathbb C[[\lambda]]$ 
be the evaluation functional at $0 \in \mathbb C^n$ and $\pi_0$ the
corresponding GNS representation on 
$\mathfrak H_0 = (\mathbb C[[\cc y^1, \ldots, \cc y^n]])[[\lambda]]$. 
Then the local and $^*$-local commutant of 
$\pi_0 (C^\infty (\mathbb C^n)[[\lambda]])$ is trivial.
\end{proposition}
\begin{proof}
We have to show that if an arbitrary $\mathbb C[[\lambda]]$-linear
endomorphism $L$ of $\mathfrak H_0$ commutes with all $\pi_0 (f)$ then 
it is a multiple of the identity. To this end we can use the
additional canonical ring structure of $\mathfrak H_0$ and in
particular the `vacuum vector' $1 \in \mathfrak H_0$. Let $L$ be such
an endomorphism. 
Since any left multiplication by elements of $\mathfrak H_0$ can be
realized as $\pi_0 (f)$ it follows that $L$ commutes with all
left multiplications and thus $L$ is itself a left multiplication by
the element $L(1) \in \mathfrak H_0$. On the other hand $L$ commutes
with  $\pi_0 (z^i) = \lambda \frac{\partial}{\partial \cc y^i}$ for
all $i$ whence $L(1)$ has to be a constant.
\end{proof}

Consider now the general case where $M$ is an arbitrary K{\"a}hler
manifold with the canonical Fedosov star product $\wick$ of Wick type
as constructed in \cite{BW97a} (see \cite{Kar96,Kar98} for another
approach to such star products of (anti-) Wick type and their
classification). In \cite[Prop.~9]{BW98a} it was shown that $\delta_p$ 
for any $p \in M$ is a positive $\mathbb C[[\lambda]]$-linear
functional for this star product. 
Now again $\supp \delta_p = \{p\}$ and the whole analysis from above
can be repeated completely analogously with the only modification that 
the Fedosov-Taylor series $\tau_p$ at $p$ enters in the analogue of
(\ref{BargmannFock}), see \cite[Thm.~5]{BW98a}. But since this map is
again surjective (which can be viewed as a sort of quantum Borel lemma
\cite[Prop.~10]{BW98a}) the same argument as above goes through. We
omit here the rather obvious details and state the final result:
\begin{theorem}
Let $(M,\wick)$ be a K{\"a}hler manifold with canonical Fedosov star
product of Wick type and let $p \in M$. Then the local 
(and $^*$-local) commutant of the GNS representation 
$\pi_p (C^\infty (M)[[\lambda]])$ induced by $\delta_p$ is trivial.
\end{theorem}

\subsection*{Schr{\"o}dinger representations on cotangent bundles}

Another important class of examples is given by cotangent bundles
whose quantization is of particular interest for physics since the
typical phase spaces are cotangent bundles of some configuration space 
manifold. In a series of papers 
\cite{BNW98a,BNW99a,BNPW98,Pfl98b,Pfl98c} 
the deformation quantization of cotangent bundles and its relation to
pseudo-differential operators and symbol calculus has been extensively
discussed and we shall investigate now the locality properties of
these star products and their representations. So let us first briefly
recall some of the basic results of \cite{BNW98a,BNW99a,BNPW98}.

One starts with a cotangent bundle $\pi: T^*Q \to Q$ over the
so-called configuration space $Q$ which can be embedded as zero
section $\iota: Q \hookrightarrow T^*Q$ in its cotangent bundle. Given 
a torsion-free connection $\nabla$ on $Q$ and a positive volume
density $\mu \in \Gamma^\infty (|\!\bigwedge^n\!|\, T^*Q)$ on $Q$ one
obtains by means of a (slightly modified) Fedosov construction firstly
the so-called standard ordered star product $\std$ which is a
homogeneous star product in the sense that the homogeneity operator 
$\mathsf H = \lambda \frac{\partial}{\partial \lambda} + \Lie_\xi$ is
a derivation of $\std$ where $\Lie_\xi$ is the Lie derivative with
respect to the Liouville vector field on $T^*Q$. Next we consider the
operator \cite[Eq.~(106)]{BNW98a}
\begin{equation}
\label{DerNeumaier}
    N = 
    \exp\left(\frac{\lambda}{2\im} 
    (\Delta + \mathsf F(\alpha))\right), 
\end{equation}
where $\Delta$ denotes the Laplacian of the semi-Riemannian metric on 
$T^*Q$ induced by the natural pairing of the horizontal and vertical
tangent spaces (locally given by 
$\Delta = \sum_k \frac{\partial^2}{\partial q^k \partial p_k}
 + \sum_{j,k,l} p_j \pi^*\Gamma^j_{kl} 
 \frac{\partial^2}{\partial p_k \partial p_l}
 + \sum_{j,k} \pi^*\Gamma^j_{jk} \frac{\partial}{\partial p_k}$, where
$\Gamma^j_{kl}$ are the Christoffel symbols of $\nabla$), and 
$\mathsf F(\alpha)$ is locally given by 
$\mathsf F(\alpha) = \sum_k \pi^*\alpha_k \frac{\partial}{\partial p_k}$
where $\alpha = \sum_k \alpha_k dq^k$ is the unique one-form such that 
$\nabla_X \mu = \alpha (X) \mu$ for $X \in \Gamma^\infty (TQ)$. 
Here we have used a canonical
(bundle) chart of $T^*Q$ but obviously the above expressions are
independent of the chart we use, see \cite{BNPW98} for a more
geometrical description of these operators.  Using $N$ as equivalence
transformation one defines the Weyl ordered product $\weyl$ by  
$f \weyl g = N^{-1} (Nf \std Ng)$, where 
$f, g \in C^\infty (T^*Q)[[\lambda]]$, generalizing thereby the
well-known Weyl-Moyal product from flat $\mathbb R^{2n}$. This star
product enjoys the following properties: firstly, we have 
$\cc{f \weyl g} = \cc g \weyl \cc f$ and, secondly, the 
$\mathbb C[[\lambda]]$-linear functional 
\begin{equation}
\label{WeylFunctDef}
    f \mapsto \omega (f) = \int_Q \iota^* f \mu,
\end{equation} 
defined on $C^\infty_0 (T^*Q)[[\lambda]]$ is positive with respect to
$\weyl$. Note that in \cite{BNW99a} this functional was defined on a larger
space namely on those functions $f$ with 
$\iota^* f \in C^\infty_0 (Q)[[\lambda]]$ but to be consistent with
our notation we shall use $C^\infty_0 (T^*Q)[[\lambda]]$ which turns
out to be still `sufficiently large'. Moreover, the Gel'fand ideal 
$\mathcal J_\omega$ is given by those functions 
$f \in C^\infty_0 (T^*Q)[[\lambda]]$ satisfying $\iota^* Nf = 0$
whence the GNS pre-Hilbert space $\mathfrak H_\omega$ is isomorphic to
the `formal wave functions' $C^\infty_0 (Q)[[\lambda]]$ on $Q$ by
\begin{equation}
\label{WeylIsoDef}
    \mathfrak H_\omega = C^\infty_0 (T^*Q)[[\lambda]] 
                         \big/ \mathcal J_\omega 
    \ni \psi_f \; \mapsto \; \iota^* Nf \in C^\infty_0 (Q)[[\lambda]].
\end{equation}
For technical reasons we choose a smooth cut-off function 
$\chi: T^*Q \to [0,1]$ such that $\chi$ is equal to $1$ in an open
neighborhood of $\iota (Q)$ and $\chi|T_q^*Q$ has compact support for 
each $q \in Q$. Then for $u \in C^\infty_0 (Q)[[\lambda]]$ the map 
$u \mapsto \psi_{\chi \pi^*u}$ is clearly the inverse of the above
isomorphism. Note that the usage of $\chi$ is necessary since we have
restricted ourselfs to $C^\infty_0 (T^*Q)[[\lambda]]$. Using these
isomorphisms  
one obtains the following explicit formula for the GNS representation
$\wrep$ of $C^\infty (T^*Q)[[\lambda]]$ on $C^\infty_0 (Q)[[\lambda]]$
\begin{equation}
\label{WeylRep}
    \wrep (f) u = \sum_{r=0}^\infty \frac{1}{r!} 
                  \left(\frac{\lambda}{\im}\right)^r
                  \sum_{i_1, \ldots, i_r} \iota^*\left(
                  \frac{\partial^r Nf}
                  {\partial p_{i_1} \cdots \partial p_{i_r}}
                  \right)
                  i_s (\partial_{q^{i_1}}) \cdots 
                  i_s (\partial_{q^{i_r}})
                  \frac{1}{r!} D^r u,
\end{equation}
where $i_s (\partial_{q^k})$ denotes the symmetric insertion of the
tangent vector $\partial_{q^k}$, and $D$ is the operator of symmetric
covariant differentiation, locally given by 
$D = \sum_k dq^k \vee \nabla_{\partial_{q^k}}$, 
see \cite[Eq.~(7)]{BNW99a}.

The first trivial observation is that the support of $\omega$ is given 
by the zero section $\iota(Q)$. Next we want to determine the support
of the equivalence class $\psi_{\chi \pi^*u}$ and expect that it
coincides with $\iota (\supp u)$. This is indeed the case as the
following simple verification shows. Let 
$0 \ne u \in C^\infty_0 (Q)[[\lambda]]$ and consider 
$\psi_{\chi \pi^* u} \in \mathfrak H_\omega$. Then for 
$g \in C^\infty_0 (T^*Q \setminus \iota (\supp u))[[\lambda]]$ we
clearly have $\omega_{\chi \pi^*u} (g) = 0$ whence 
$\supp \psi_{\chi\pi^* u} \subseteq \iota (\supp u)$. The converse
inclusion is also true. Let $O \subseteq T^*Q$ be open such that 
$O \cap \iota (\supp u) \ne 0$ then $\pi (O \cap \iota (\supp u))$ is
open and non-empty in $Q$. Choose a non-negative function 
$0 \ne v \in C^\infty_0 (\pi (O \cap \iota (\supp u)))[[\lambda]]$
then clearly $\omega_{\chi \pi^* u} (\tilde\chi \pi^*v) > 0$ where
$\tilde \chi$ is a suitable smooth bump function equal to $1$ in a
neighborhood of $O \cap \iota (\supp v)$ such that 
$\supp (\tilde\chi \pi^* v) \subseteq O$. Thus 
$\iota(\supp u) \subseteq \psi_{\chi \pi^*u}$ whence we have shown the
following lemma:
\begin{lemma}
Let $\psi_f \in \mathfrak H_\omega$ then 
$\supp \psi_f = \iota (\supp (\iota^* Nf))$.
\end{lemma}
Thus the abstract definition of $\supp \psi_f$ coincides with the
usual geometric support of the corresponding formal wave function on
$Q$ (embedded in $T^*Q$). Hence a local operator on 
$\mathfrak H_\omega$ corresponds under the isomorphism
(\ref{WeylIsoDef}) to a local operator on $C^\infty_0 (Q)[[\lambda]]$
in the usual sense. Thus the general statement of
Theorem~\ref{LocalGNSOp} that the GNS representation automatically
yields local operators is 
manifested  here by the fact that clearly $\wrep (f)$ is a formal
series of differential operators and thus local. This observation
enables us to compute the local commutant of 
$\wrep (C^\infty (T^*Q)[[\lambda]])$.
\begin{theorem}
\label{SchroedingerTheo}
The local commutant $(\wrep (C^\infty (T^*Q)[[\lambda]]))'_{\Loc}$ of
the Schr{\"o}dinger representation is trivial if and only if $Q$ is
connected. In general $(\wrep (C^\infty (T^*Q)[[\lambda]]))'_{\Loc}$
is isomorphic to $H_0 (Q)[[\lambda]]$.
\end{theorem}
\begin{proof}
Let 
$A: C^\infty_0 (Q)[[\lambda]] \to C^\infty_0 (Q)[[\lambda]]$ be a
local $\mathbb C[[\lambda]]$-linear operator commuting with all
$\wrep(f)$. Since $A$ is $\mathbb C[[\lambda]]$-linear it is of the
form $A = \sum_{r=0}^\infty \lambda^r A_r$ and clearly all operators
$A_r: C^\infty_0 (Q) \to C^\infty_0 (Q)$ are local again. Since $A$
commutes with all left multiplications by functions 
$u \in C^\infty_0 (Q)$ the lowest order $A_0$ commutes with all such
left multiplication. On the other hand by Petree's theorem 
(see e.g.\ \cite[p.~176]{KMS93}) the locality of $A_0$ implies that
around any point $q \in Q$ there is a chart such that $A_0$
restricted to this chart is a \emph{differential} operator. Putting
this together we see that $A_0$ has to be even of order zero in this
and hence in any chart, i.e.\ $A_0$ is a left multiplication by a
function $a_0 \in C^\infty (Q)$ itself. On the other hand $A_0$
commutes with any Lie derivative $\Lie_X$ which can be obtained by 
$\wrep (\hat X)$ minus some left multiplications where $\hat X$ is the 
function linear in the momentum variables given by 
$\hat X (\alpha_q) = \alpha_q (X_q)$ where $\alpha_q \in T^*_q Q$ is a 
point in $T^*Q$ and $X \in \Gamma^\infty (TQ)$ is a vector field.
Thus $a_0$ has to be constant on each connected component. 
Induction on $r$ completes the proof.
\end{proof}
\begin{remark}
\begin{enumerate}
\item In \cite{BNPW98,BNW98a,BNW99a} several generalizations for
      $\weyl$ and $\omega$ have been made: firstly, one can associate
      to any projectable Lagrangian submanifold $L$ of $T^*Q$ a
      functional $\omega_L$ which in the corresponding GNS
      representation induces the WKB expansion of a Hamiltonian $H$
      satisfying the Hamilton-Jacobi equation $H|L = E$ for some
      energy value $E$, see e.g.\ \cite{BaWei95}. 
      In this case the support of the functional
      $\omega_L$ is $L$ and the GNS pre-Hilbert space is isomorphic to 
      $C^\infty_0 (L)[[\lambda]]$ such that the abstract support again
      corresponds to the geometric support on the Lagrangian
      sub-manifold. Secondly one can also incorporate a `magnetic
      field' as an additional closed two-form on $Q$ pulled back to
      $T^*Q$ and added to the canonical symplectic form. The
      corresponding star products have GNS representations on a
      Hermitian line bundle over $Q$ in case where the magnetic field
      satisfies an additional integrality condition. Then the GNS
      pre-Hilbert space is isometric to the sections of this line
      bundle with compact support and the locality structure of the
      abstract quotient coincides under this isomorphism with the
      usual notion of the support of the sections. Finally in both
      cases the local commutant is again trivial if and only if $Q$ is 
      connected. We shall not carry out this in detail since the proof
      works completely analogously.
\item The above examples show that these various kinds of
      Schr{\"o}dinger representations as well as the Bargmann-Fock
      representation with their local commutants indeed behave
      completely different to the thermal KMS representations. In
      particular in the case of the Schr{\"o}dinger representation the 
      connected components of the configuration space $Q$ behave like
      super selection rules and in fact are the only ones, see also
      Corollary~\ref{SuperSelectionCor}.
\item Finally this example shows that the notion of irreducible
      representations does not seem to be appropriate for deformation
      quantization since clearly
      $\mathfrak H_\omega (T^*O) \cong \mathcal A_0 (O)$ is an
      invariant subspace for all open $O \subseteq Q$. Thus the
      characterization by local commutants is more suitable. Moreover,
      the locality structure of the pre-Hilbert space 
      $\mathfrak H_\omega (Q)$, i.e.\ the fact that 
      $\mathfrak H_\omega (O) \bot \mathfrak H_\omega (O')$ for 
      $O \cap O' = \emptyset$, can be understood as a
      \emph{consequence} of the locality structure of the observable
      algebra. If one heuristically thinks of formal deformation
      quantization as an `asymptotic expansion' of some convergent
      theory, it seems that having a (non-)trivial commutant is a more
      `rigid' property of a representation with respect to asymtotic
      behaviour than being (ir-)reducible. 
\end{enumerate} 
\end{remark}

%
%

\section{Strong topologies and von Neumann algebras}
\label{NeuSec}

In this section we shall investigate further similarities between
local operators in formal GNS representations and bounded operators on
complex Hilbert space and end up with certain analogues of von
Neumann's double commutant theorem.

Let us first introduce the notion of an \emph{approximate identity}
borrowed from $C^*$-algebra theory \cite[Def.~2.2.17.]{BR87}. Let
$\{O_n\}_{n \in \mathbb N}$ be a sequence of open subsets of $M$ such
that each $O_n$ has compact closure $O_n^{\mathrm {cl}}$ contained in
$O_{n+1}$ and such that $M = \bigcup_{n\in \mathbb N} O_n$. Furthermore
let $\chi_n \in C^\infty_0 (O_{n+1})$ be a smooth function such that
$\chi_n|O_n^{\mathrm {cl}} = 1$ for all $n \in \mathbb N$. 
Then $(O_n, \chi_n)_{n\in \mathbb N}$ is called an approximate
identity. Note that there always exists such an approximate identity
and in the case where $M$ is compact we simply may choose $\chi_n = 1$
and $O_n = M$ for all $n \in \mathbb N$.

Now let $(M, *)$ be a Poisson manifold with local star product. Then
we consider the space $C^\infty_0 (M)[[\lambda]]$ and its
local $\mathbb C[[\lambda]]$-linear endomorphisms.
Let $\AL$ and $\AR$ denote all those endomorphisms
obtained by left and right multiplication with elements of 
$C^\infty (M)[[\lambda]]$, respectively. Clearly left and right
multiplications commute whence $\AL \subseteq \AR'$ and 
$\AR \subseteq \AL'$, where we may take the commutant in all 
$\mathbb C[[\lambda]]$-linear endomorphisms or in the local ones. Note 
that canonically $\AL \cong (C^\infty (M)[[\lambda]], *)$ and 
$\AR \cong (C^\infty (M) [[\lambda]], *)^{\mathrm {op}}$,
respectively.
\begin{proposition}
\label{ALARCommutantLem}
Let $(M, *)$ be a Poisson manifold with local star product then
\begin{equation}
    \AL' = \AR 
    \quad
    \textrm { and }
    \quad
    \AR' = \AL
\end{equation}
\end{proposition}
The proof is trivial if $M$ is compact since in this case 
$1 \in C^\infty_0 (M)[[\lambda]]$ and in the non-compact case one uses 
an approximate identity.

Though the centers of $(C^\infty (M)[[\lambda]], *)$, $\AL$, and
$\AR$, respectively, could be rather large in the general Poisson case 
the centers are known to be trivial in the case where $M$ is a
connected symplectic manifold.
\begin{lemma}
\label{TrivialCenterLem}
Let $(M, *)$ be a connected symplectic manifold with local star
product. Then the centers of $(C^\infty (M)[[\lambda]], *)$, $\AL$, and
$\AR$ are trivial.
\end{lemma}

In (\ref{KMSAdjoint}) we have already noticed that for a KMS functional the 
adjoints of all right multiplications exits with respect to the 
induced Hermitian product on 
$\mathfrak H_\KMS \cong C^\infty_0 (M)[[\lambda]]$.
As a corollary we obtain that both $\AL$ and $\AR$, viewed as
subalgebras of $\BLoc (\mathfrak H_\KMS)$, are factors if $M$ is
connected: 
\begin{corollary}
Let $(M, *)$ be a connected symplectic manifold with local star
product and let $\omega_\KMS$ be a positive KMS functional as in
(\ref{KMSFunctional}). Then $\AL$ and $\AR$ are factors in 
$\BLoc(\mathfrak H_\KMS)$.
\end{corollary}

Since as $\mathbb C[[\lambda]]$-modules canonically 
$\mathfrak H_\omega \cong C^\infty_0 (M)[[\lambda]]$ (including the
corresponding locality structures) when $\omega$ is faithful, we shall 
now investigate $C^\infty_0 (M)[[\lambda]]$ and the local operators of 
this space in more detail. In particular we are interested in
topological properties of the local operators 
$\Loc (C^\infty_0 (M)[[\lambda]])$ and the relation of the topological 
closures with double commutant closures. To this end we have to
specify the topologies we want to use, but firstly it will be necessary 
to enlarge the framework to more general series. We need (at
least) formal Laurent series in $\lambda$, see App.~\ref{FormalApp}
for definitions. It is clear that all definitions and results are also 
valid in this setting \emph{if} we require not only 
$\LS{\mathbb C}$-linearity but in addition also $\lambda$-adic
continuity of all involved maps, as e.g.\ the positive functional
$\omega$, the local operators etc. This is crucial in view of 
Lemma~\ref{LambdaContLem}, and the possible complications cannot be
seen in the framework of formal power series since here 
$\mathbb C[[\lambda]]$-linearity implies $\lambda$-adic 
continuity, see e.g.\ \cite[Prop.~2.1]{DL88}. Thus let 
$\omega: \LS{C^\infty_0 (M)} \to \LS{\mathbb C}$ from now on be a
$\lambda$-adically continuous, positive and $\LS{\mathbb C}$-linear
functional and denote by $\Loc(\mathfrak H_\omega)$ those
$\LS{\mathbb C}$-linear endomorphisms of $\mathfrak H_\omega$ which
are local and $\lambda$-adically continuous. Then it is clear due to
Lemma~\ref{LambdaContLem} that $\Loc\left(\LS{C^\infty_0 (M)}\right)$ 
is given by $\LS{\Loc (C^\infty_0 (M))}$ where we denote by 
$\Loc (C^\infty_0 (M))$ the (usual) $\mathbb C$-linear local operators  
on $C^\infty_0(M)$. The $\lambda$-adic topology of 
$\LS{\Loc (C^\infty_0(M))}$ will be somehow too fine for the 
study of von Neumann algebras (similar to the 
norm-topology of bounded operators on complex Hilbert spaces) and thus 
we are using the following \emph{strong operator topology} analogously 
to the usual situation in complex Hilbert spaces. A basis 
of open neighborhoods of $0 \in \LS{\Loc (C^\infty_0(M))}$ will be
given by
\begin{equation}
    O_{f_1, \ldots, f_n; \epsilon} 
    := \left\{ A \in \LS{\Loc(C^\infty_0(M))} \: | \:
    \forall l = 1, \ldots, k: \;
    \varphi (A (f_l)) < \epsilon \right\},
\end{equation}
where $\epsilon > 0$, $f_1, \ldots, f_k \in \LS{C^\infty_0 (M)}$, and
$\varphi$ is the $\lambda$-adic absolute value.
Then a sequence $A_n$ of local operators converges \emph{strongly}, 
i.e.\ with respect to this topology, to $A$ if and only if for all
$f \in \LS{C^\infty_0 (M)}$
\begin{equation}
    A_n f \to A f
\end{equation}
in the $\lambda$-adic topology of $\LS{C^\infty_0 (M)}$. Clearly, if
$A_n \to A$ in the $\lambda$-adic topology then also $A_n \to A$ in
the strong topology. The following example shows that the converse is
not true in general:
\begin{example}
\label{DiffLocEx}
Let $M = \mathbb R$ and let $\chi_0$ be a smooth function having
support in $[0,1]$. Define $\chi_n (x) := \chi_0 (x-n)$ for 
$n \in  \mathbb N$ and let $A_n$ be the left multiplication by 
$\chi_n$. It follows that $A_n$ converges strongly to $0$ 
but it does not converge in the $\lambda$-adic topology.
\end{example}
From the general statement in Prop.~\ref{CauchyLocProp} we see that
$\LS{\Loc (C^\infty_0 (M))}$ is complete with respect to this
topology. Hence it makes sense to ask whether a double commutant
coincides with a topological closure in order to find at least for
particular cases an analogue to von Neumann's double commutant
theorem.

We shall now turn again to faithful functionals $\omega$ since in this
case $\mathfrak H_\omega \cong \LS{C^\infty_0 (M)}$. Thus it will be
sufficient to consider the latter space. Moreover, we
consider the algebra $\ALR$ which is generated by all left and right
multiplications $\AL$ and $\AR$. Note that we have a canonical
surjective morphism $\AL \otimes_{\LS{\mathbb C}} \AR \to \ALR$,
simply given by 
$\mathsf L_f \otimes \mathsf R_g \mapsto \mathsf L_f \mathsf R_g$ which
is \emph{not} injective since whenever 
$\supp f \cap \supp g = \emptyset$ we have 
$\mathsf L_f \mathsf R_g = 0$ but in general
$\mathsf L_f \otimes \mathsf R_g \ne 0$. Moreover, $\AL$ and $\AR$ are
canonically embedded in both 
$\AL \otimes_{\LS{\mathbb C}} \AR$ and $\ALR$.

In the general Poisson case the local commutant of $\ALR$ can be
rather big (since its center can be rather big) but for the connected
symplectic case the commutant is trivial due to 
Lem.~\ref{TrivialCenterLem} and Prop.~\ref{ALARCommutantLem}.
\begin{lemma}
Let $(M, *)$ be a connected symplectic manifold with local star
product. Then $(\ALR)'_\Loc = \LS{\mathbb C} \id$ whence 
$(\ALR)''_\Loc = \LS{\Loc (C^\infty_0 (M))}$.
\end{lemma}
Thus we may now ask whether $\ALR$ is dense in its
double commutant. With the above strong topology this is indeed the
case:
\begin{theorem}
\label{vonNeumannTheo}
Let $(M,*)$ be a connected symplectic manifold with local star product.
Then the completion of $\ALR$ in the
strong operator topology is $\LS{\Loc(C^\infty_0 (M))}$.
\end{theorem}
\begin{proof}
For a given local operator $L \in \Loc(\LS{C^\infty_0(M)})$ 
we have to construct a sequence of elements $A_n \in \ALR$
such that $A_n \to L$ strongly.
First consider a formal series of local operators $D$ whose
coefficients have support in one common compact set contained in
some open subset $O$ of $M$. By Peetre's theorem we know, that each
coefficient is then a differential operator and by 
$\LS{\mathbb C}$-linearity we may assume $o(D) = 0$ whence we write 
$D = \sum_{r=0}^\infty \lambda^r D_r$. Consider now 
$g \in C^\infty_0 (O)$ then $\ad (g) = \mathsf L_g - \mathsf R_g$ has
order $\ge 1$ and starts in lowest order with the Poisson bracket
$\ad (g) = \im\lambda \{g, \cdot \} + \ldots$. 
Since the Poisson bracket is non-degenerate in the
\emph{symplectic} case, we can obtain by suitable choice of finite 
algebraic combinations of left multiplications and commutators with
elements in $C^\infty_0 (O)$ any differential operator $D_0$ up to
higher orders in $\lambda$, if we allow for division by
finitely many powers of $\lambda$ (actually by $\lambda^k$ if the
order of differentiation is $k$). Thus we obtain $A_0 \in \ALR$
having support in $O$ such that $o (D - A_0) \ge 1$ and by induction
we find for any $n \in \mathbb N$ an element $A_n \in \ALR$ having
support in $O$ such that $o(D-A_n) \ge n$ since the higher orders of
the operators of left and right multiplications with elements having
compact support are differential operators due to the locality of $*$
and Peetre's theorem.

Now let $L$ be a local operator and let 
$(O_n, \chi_n)_{n \in \mathbb N}$ be an approximate identity. Then
$\chi_n L$ is still a local operator having compact support in
$O_{n+1}$. Thus we can find a sequence $A_n \in \ALR$ such that 
$o(\chi_n L - A_n) \ge n+1$ for all $n \in \mathbb N$. We claim 
$A_n \to L$ strongly. To prove this let $f \in \LS{C^\infty_0(M)}$
where we can assume by $\LS{\mathbb C}$-linearity that 
$o(f) = 0$. Write $f = \sum_{r=0}^\infty \lambda^r f_r$ then for any
$k \in \mathbb N$ there is a $N \in \mathbb N$ such that $N > k$ and 
$\supp f_1 \cup \cdots \cup \supp f_k \subseteq O_n$ for all $n \ge N$.
Thus $\chi_n L f_r = L f_r$ for $r=1, \ldots, k$ whence 
$o(\chi_n L f - Lf) \ge k$ for all $n \ge N$. On the other hand 
$o(\chi_n L - A_n) \ge n+1$ uniformly whence also 
$o (\chi_n L f - A_n f) \ge n+1$ since $o(f) = 0$ for all $n \ge N$.  
Then the strong triangle inequality for the order implies that 
$o(Lf - A_n f) \ge k$ for all $n \ge N$ whence indeed $A_n \to L$ in
the strong operator topology.
\end{proof}

In order to get the full analogy of von Neumann's double commutant
theorem we have to take into account the $^*$-involution too. Hence we
define the $^*$-strong operator topology by specifying the following
basis of open neighborhoods of $0$:
\begin{equation}
    O_{f_1, \ldots, f_k; \epsilon} :=
    \left\{
    A \in \BLoc (\LS{C^\infty_0 (M)}) \; | \; 
    \forall l = 1, \ldots, k: \;
    \varphi (A(f_l)) < \epsilon
    \textrm { and }
    \varphi (A^*(f_l)) < \epsilon \right\},
\end{equation}
where $\epsilon > 0$ and 
$f_1, \ldots, f_k \in \LS{C^\infty_0 (M)}$, $k \in \mathbb N$. 
Then $A_n$ converges $^*$-strongly to $A$ if and only if the sequences
$A_n f$ and $A_n^* f$ converge to $Af$ and $A^*f$, respectively,
in the $\lambda$-adic topology 
for all $f \in \LS{C^\infty_0 (M)}$. Note that this topology
incorporates now the GNS Hermitian product of the faithful functional
$\omega$. Clearly the $^*$-strong operator topology is finer than the
strong operator topology whence $A_n \to A$ $^*$-strongly implies 
$A_n \to A$ strongly but the reverse needs not to be true. 
\begin{proposition}
Let $(M, *)$ be a Poisson manifold with local star product and 
$\omega: \LS{C^\infty_0 (M)} \to \LS{\mathbb C}$ be a faithful,
positive, $\lambda$-adically continuous, and $\LS{\mathbb C}$-linear
functional. Then the space $\BLoc (\LS{C^\infty_0 (M)})$ is
complete in the $^*$-strong operator topology.
\end{proposition}
\begin{proof}
Let $A_n \in \BLoc (\LS{C^\infty_0 (M)})$ be a $^*$-strong Cauchy
sequence. Since $A_n \in \Loc (\LS{C^\infty_0 (M)})$ this implies that
$A_n$ is Cauchy with respect to the strong operator topology, too, and
by Prop.~\ref{CauchyLocProp} convergent to some
$A \in \Loc (\LS{C^\infty_0 (M)})$. Similarly $A_n^*$ converges
strongly to some $B \in \Loc (\LS{C^\infty_0 (M)})$. Thus it remains
to show that $B = A^*$ and $A_n \to A$ $^*$-strongly. But this is a
simple verification using the $\lambda$-adic continuity of $\omega$
and the Hermitian product. 
\end{proof}

We conclude from Theorem~\ref{vonNeumannTheo} and this proposition
that for the connected symplectic case with a KMS functional the 
$^*$-strong completion of $\ALR$ is $\BLoc (\LS{C^\infty_0 (M)})$:
\begin{corollary}
Let $(M, *)$ be a connected symplectic manifold and $\omega_\KMS$ a
positive KMS functional as in (\ref{KMSFunctional}). Then the
$^*$-strong completion of $\ALR$ is $\BLoc (\LS{C^\infty_0 (M)})$.
\end{corollary}

Note that for the above proofs both the usage of formal Laurent series
and the non-degeneracy of the Poisson bracket were crucial. It remains
an open and interesting problem whether and how the above
theorem can be extended to the general Poisson case. Here a possible
degeneracy of the Poisson bracket in certain directions may be
compensated by higher orders of the star product. As an example one
can consider a symplectic manifold $M$ with star product $*$. Then the
substitution $\lambda \mapsto \lambda^2$ provides a star product for
the Poisson bracket which vanishes identically but clearly the above
theorem is still valid in this case. 
For the general case the above strong operator topology may still be
too fine and in order to find a coarser topology one might have 
to take into account the locally convex topology of $C^\infty_0(M)$
too.

An analogous theorem is valid for the Schr\"odinger-like GNS
representations on cotangent bundles. With the notation from
Section~\ref{ExSec} and the extension to formal Laurent series we have
$\mathfrak H_\omega \cong \LS{C^\infty_0 (Q)}$ and the representation
is given by (\ref{WeylRep}). Since here (if $Q$ is connected) the
local commutant of $\wrep (\LS{C^\infty_0(T^*Q)})$ is already trivial,
we expect that the strong closure yields all local operators on
$\mathfrak H_\omega$. This is indeed the case:
\begin{theorem}
Let $Q$ be a connected manifold endowed with a torsion-free connection
and a positive density. Then for the corresponding Weyl ordered star
product algebra and its Schr\"odinger representation as in
Section~\ref{ExSec} the local operators $\Loc (\LS{C^\infty_0 (Q)})$ are
the completion of $\wrep (\LS{C^\infty (T^*Q)})$ in the strong
operator topology.
\end{theorem}
\begin{proof}
Since we allowed for finitely many negative powers of $\lambda$ we
notice from (\ref{WeylRep}) that any differential operator on
$C^\infty_0 (Q)$ can be expressed as $\wrep (f)$ with a suitably
chosen $f \in \LS{C^\infty (T^*Q)}$. Thus the claim easily follows
from Example~\ref{CauchyDiffLocEx} and Prop.~\ref{CauchyLocProp}.
\end{proof}

We conclude this section with a few remarks: Due to the particular and
simple form of the Hermitian product in the Schr\"odinger
representation one observes that in this case any local operator has
an adjoint. In the general case (even for faithful positive
functionals) this is not obvious, whence in this case the $^*$-strong
operator topology is needed. Moreover, if we enlarge the 
framework to formal CNP series
then the strong operator topology can also be written by use of a norm
topology of the underlaying Hilbert space over $\CNP{\mathbb C}$ since
in this case we can define a $\CNP{\mathbb C}$-valued norm of $\psi_f$
by $\| \psi_f \| := \sqrt{\SP{\psi_f, \psi_f}}$. Many aspects of
such Hilbert spaces over the field $\CNP{\mathbb C}$ were discussed in
\cite{BW98a}. One aim to do this could be a `formal spectral theory'
within the local operators in order to compute formal spectra and
compare them with asymptotic expansion of their convergent
counterparts (if there exists a convergent counterpart).

%
%

\section{Tomita-Takesaki theory}
\label{ToTaSec}

Since the concept of KMS functionals can be formulated for deformation
quantization the question for an analogue of the usual Tomita-Takesaki
theory arises naturally. It turns out that the deformed algebras allow
indeed for such an analogue which will be surprisingly simple: it can
be formulated purely algebraically and the usual functional-analytical
difficulties do not occur.

Let $(M, *)$ be a connected symplectic manifold and let 
$H = \sum_{r=0}^\infty \lambda^r H_r \in C^\infty (M)[[\lambda]]$ be a
real Hamiltonian and $\beta \in \mathbb R$ an inverse
temperature. Then we denote the corresponding KMS functional by 
$\omega_\KMS$ which is given as in (\ref{KMSFunctional}). The GNS
pre-Hilbert space $\mathfrak H_\KMS$ is then isomorphic to 
$C^\infty_0 (M)[[\lambda]]$ and in this section we shall always
identify them. Using the same notation as for the usual
Tomita-Takesaki theory, see e.g.\ \cite[Sect.~2.5.2]{BR87}, we define
the $\mathbb C[[\lambda]]$-anti-linear operator 
$S: C^\infty_0 (M)[[\lambda]] \to C^\infty_0 (M)[[\lambda]]$ by
$Sf := \cc f$. Since the space $C^\infty_0 (M)[[\lambda]]$ is already
complete with respect to the $\lambda$-adic topology $S$ is defined on
the whole GNS representation space which drastically simplifies the
approach. By a simple computation we see that the operator 
$F: f \mapsto \Exp (-\beta H) * \cc f * \Exp (\beta H)$ is the unique
$\mathbb C[[\lambda]]$-anti-linear adjoint of $S$, i.e.\ we have
$\SP{f, Sg}_\KMS = \cc{\SP{Ff, g}}_\KMS$ for all 
$f, g \in C^\infty_0 (M)[[\lambda]]$. Thus we define the 
$\mathbb C[[\lambda]]$-linear operator $\Delta := FS$ as usual and
obtain the explicit expression
\begin{equation}
    \Delta f = \Exp (-\beta H) * f * \Exp (\beta H)
             = \mathsf L_{\Exp (-\beta H)}
               \mathsf R_{\Exp (\beta H)} f.
\end{equation}
Clearly $\Delta$ is positive with respect to the KMS Hermitian
product, i.e.\ we have $\SP{f, \Delta f} \ge 0$ by a simple
computation. Moreover, $\Delta$ is obviously invertible with inverse
$\Delta^{-1} = \mathsf L_{\Exp (\beta H)}\mathsf R_{\Exp (-\beta H)}$,
and for all $z \in \mathbb C[[\lambda]]$ we define
$\Delta^z := \mathsf L_{\Exp (-z\beta H)}\mathsf R_{\Exp (z\beta H)}$ 
such that $\Delta^z\Delta^{z'} = \Delta^{z+z'}$ and 
$\Delta^0 = \id$. Hence we can define the following 
$\mathbb C[[\lambda]]$-anti-linear operator 
$J: C^\infty_0 (M)[[\lambda]] \to C^\infty_0 (M)[[\lambda]]$ by
\begin{equation}
    J := S \Delta^{-\frac{1}{2}}.
\end{equation}
A straightforward computation yields the explicit expression 
\begin{equation}
    Jf = \Exp \left(-\frac{\beta}{2} H\right) * 
         \cc f * 
         \Exp \left(\frac{\beta}{2} H\right), 
\end{equation}
whence in particular $J^2 = \id$ as well as 
$\SP{Jf, Jg} = \SP{g, f}$. Hence $J$ is anti-unitary with
$J = J^* = J^{-1}$. Finally observe that 
$J \Delta^{\frac{1}{2}} J = \Delta^{-\frac{1}{2}}$ and 
$S^2 = \id = F^2$.
Analogously to the usual Tomita-Takesaki theory we call $J$ the
\emph{modular conjugation} and $\Delta$ the \emph{modular operator}
associated to the KMS functional $\omega_\KMS$. It remains to give a
reasonable definition of the \emph{modular group}, i.e.\ the 
`time development' induced by the modular operator. Thus we have to
find a reasonable definition for $\Delta^{\frac{\im t}{\beta\lambda}}$
where the obvious problem comes from the $\lambda$ in the
denominator. Note that this is necessary to get the physical
dimensions right. Thus a naive definition as for $\Delta^z$ is not
possible since the star exponential 
$\Exp \left(\frac{\im t}{\beta\lambda} H \right)$ would not be 
well-defined in the category of formal power series. Since we do not
want to leave this framework we have to give an alternative
definition. To motivate this, let us proceed heuristically for a
moment: if $\Delta^{\frac{\im t}{\beta\lambda}}$ were a one-parameter
group of unitaries of $\mathfrak H_\KMS$ we could try to differentiate
it in order to find a differential equation which we can afterwards
solve to define the modular group. Since obviously a formal logarithm
of $\Delta$ is defined, namely $\ln \Delta = -\beta \ad (H)$, we get
the following differential equation (in a strong sense, i.e.\ after
applying to some vector $f \in C^\infty_0 (M)[[\lambda]]$):
\begin{equation}
   \frac{d}{dt} \Delta^{\frac{\im t}{\beta\lambda}} f
   =
   - \frac{\im}{\lambda} \ad(H) \Delta^{\frac{\im t}{\beta\lambda}} f
\end{equation}
But this equation now makes perfectly sense since the operator 
$\ad (H)$ is of order $\lambda$ cancelling the $\lambda$ in the
denominator. Moreover, this equation, viewed as an equation for a
time-dependent map $U_t = \Delta^{\frac{\im t}{\beta\lambda}}$ has
indeed a solution, namely the one-parameter group of automorphisms
$A_{-t}$ induced by the Heisenberg equation corresponding to the
Hamiltonian $H$ \emph{provided} the flow of the classical Hamiltonian
vector field exists for all times $t \in \mathbb R$, 
see e.g.\ \cite[App.~B]{BNW99a}. Thus we are led to the
following definition: Assume in addition that the classical flow of
the Hamiltonian vector field of $H_0$ exists (in this case the KMS
functional is also called a dynamical KMS functional), then we
\emph{define} the 
modular group by $U_t := A_{-t}$, where $A_t$ is the one-parameter
group of automorphisms of the quantum time development. Here the minus
sign is due to the fact that $f \in C^\infty_0 (M)[[\lambda]]$ is now
to be considered as state and not as observable. Using the properties
of $A_t$ one proves by direct computation that $U_t$ is indeed a
one-parameter group of unitaries of the GNS pre-Hilbert space
$C^\infty_0 (M)[[\lambda]]$, i.e.\ we have 
$\SP{U_t f, U_t g}_\KMS = \SP{f, g}_\KMS$ for all
$f, g \in C^\infty_0 (M)[[\lambda]] $ and all $t \in \mathbb R$.
Together with the fact that $\AL' = \AR$
(Prop.~\ref{ALARCommutantLem}) we can now formulate the
analogue of the Tomita-Takesaki theorem:
\begin{theorem}
With the notation from above we have:
\begin{enumerate}
\item The $\mathbb C[[\lambda]]$-anti-linear map 
      \begin{equation}
          \AL \ni \mathsf L_f  \; \mapsto \;
          J \mathsf L_f J =
          \mathsf R_{\Exp\left(-\frac{\beta}{2} H \right)
          * \cc f * \Exp\left(-\frac{\beta}{2} H \right)}
          \in \AL' = \AR
      \end{equation}
      is a bijection whence $J \AL J = \AL'$.
\item For all $z \in \mathbb C[[\lambda]]$ one has
      \begin{equation}
          \Delta^z \AL \Delta^{-z} = \AL.
      \end{equation}
\item If in addition the flow of the Hamiltonian vector field of $H_0$
      exists for all times $t \in \mathbb R$, 
      i.e.\ the KMS functional is a dynamical KMS functional, 
      then one has for all $t \in \mathbb R$
      \begin{equation}
          U_t \AL U_{-t} = \AL.
      \end{equation}
\end{enumerate}
\end{theorem}
\begin{proof}
With the above definitions the proof is a simple computation.
\end{proof}

This surprisingly simple and algebraic proof (except of the definition
of $U_t$) suggests once more that the algebras of deformation
quantization of finite-dimensional symplectic manifolds correspond
heuristically to the most simple counterpart in the usual theory of
von Neumann algebras, i.e.\ to the type I.
From the physical point of view this can be understood since 
we have only dealt with finitely many degrees of freedom, whence 
a type I$_\infty$ or even I$_n$ is expected for the quantum mechanical
description. Hence the
full complexity of the usual Tomita-Takesaki theory is not yet reached
and hence it would be of major interest to find formulations for
infinitely many degrees of freedom, where both either a quantum field
theoretical or a thermodynamical approach would be very interesting.

On the other hand, the above formulation deals only with the
symplectic case. For the general Poisson case many of the above
results were not true in general or are rather non-obvious as e.g.\ 
the existence or uniqueness of traces. In \cite{Wei97} a classical
version for Poisson manifolds is discussed and it would be very
interesting to find analogues to these and the above statements for
the quantized versions, too.

%
%

\appendix

%
%

\section{Pre-Hilbert spaces over ordered rings and the GNS construction}
\label{GNSApp}

For the reader's convenience we shall summarize here some facts on
pre-Hilbert spaces over ordered rings and related GNS constructions.
See \cite{BW98a,BNW99a} for a detailed exposition and proofs. For the
well-known case of $C^*$-algebras over complex numbers, see
e.g.\ \cite{BR81,BR87,Con94,Haa93} and also \cite{Schmue90} for more
general complex $^*$-algebras.

First recall that a commutative associative ring $\field R$ with 
$1 \ne 0$ is called \emph{ordered} with positive elements 
$\field P \subset \field R$ if $\field R$ is the disjoint union 
$\field R = -\field P \mathop{\dot\cup} \{0\} \mathop{\dot\cup} \field P$
and $\field P$ is closed under addition and multiplication. If 
$\field R$ is ordered  then it is of characteristic zero, i.e.\  
$n1 = 1 + \cdots + 1 \ne 0$ for all $n \in \mathbb Z$ and it has no
zero divisors. The quotient field $\hat{\mathsf R}$ of $\mathsf R$
becomes an ordered field such that the usual embedding of $\mathsf R$
in $\hat{\mathsf R}$ is compatible with the ordering. 
Now let $\field R$ be an ordered ring then we consider 
$\field C = \field R \oplus \im \field R$ where we endow 
$\field C$ with a ring structure by requiring $\im^2 = -1$. Then 
$\field C$ is again an associative commutative ring with 
$1 \ne 0$ and has no zero divisors. 
Elements in $\field C$ are written as $z = a + \im b$, where 
$a,b \in \field R$, and $\field R$ is embedded in $\field C$ via
$a \mapsto a + \im 0$. Complex conjugation in $\field C$ is defined as
usual by $z = a + \im b \mapsto \cc z = a - \im b$. Then 
$z \in \field C$ is an element in $\field R$ if and only if 
$z = \cc z$ and clearly $\cc z z \ge 0$ and $\cc zz = 0$ if and only
if $z =0$.

A \emph{pre-Hilbert space} over such a ring $\field C$ is defined to
be a $\field C$-module $\mathfrak H$ endowed with a $\field C$-valued
Hermitian product, i.e.\ a map 
$\SP{\cdot,\cdot} : \mathfrak H \times \mathfrak H \to \field C$ 
satisfying the following axioms: $\SP{\cdot,\cdot}$ is $\field C$-linear 
in the second argument, $\SP {\phi,\psi} = \cc {\SP{\psi,\phi}}$ for all 
$\psi, \phi \in \mathfrak H$, and $\SP{\cdot,\cdot}$ is positive,
i.e.\ 
$\SP{\psi, \psi} \ge 0$ and $\SP{\psi,\psi} = 0$ implies $\psi =0$ for
all $\psi \in \mathfrak H$. Then the Hermitian product satisfies the 
\emph{Cauchy-Schwarz inequality}
\begin{equation}
\label{CauchySchwarzI}
    \SP{\phi,\psi}\cc{\SP{\phi,\psi}} \le \SP{\phi,\phi}\SP{\psi,\psi},
    \qquad
    \phi,\psi \in \mathfrak H.
\end{equation}
If $\mathfrak H, \mathfrak K$ are pre-Hilbert spaces over $\field C$
then a $\field C$-linear map $U: \mathfrak H \to \mathfrak K$ is
called \emph{isometric} if $\SP{U\phi, U\psi} = \SP{\phi,\psi}$ for
all $\phi,\psi \in \mathfrak H$ and \emph{unitary} if in addition $U$
is surjective. Note that if $U$ is isometric then it is in particular
injective, whence a unitary map is invertible and the inverse of a
unitary map is unitary again.

Next we consider the possibility to define adjoints of
endomorphisms. Let $A: \mathfrak H \to \mathfrak H$ be a
$\field C$-linear endomorphism of  
a pre-Hilbert space $\mathfrak H$ over $\field C$. Then a 
$\field C$-linear endomorphism $B: \mathfrak H \to \mathfrak H$ is
called \emph{adjoint} of $A$, written as $B = A^*$, if for all 
$\phi, \psi \in \mathfrak H$
\begin{equation}
\label{AdjointRelation}
    \SP{B\phi, \psi} = \SP{\phi, A\psi}.
\end{equation}
In general the existence of such adjoints is far from being obvious
but if $A$ has an adjoint then it is unique. Moreover, if $A^*$ and
$B^*$ exist then $(aA + bB)^*$, $(AB)^*$, and $(A^*)^*$ exist and are
given by 
\begin{equation}
\label{Adjoints}
    (aA+bB)^* = \cc a A^* + \cc b B^*,
    \qquad
    (AB)^* = B^* A^*,
    \qquad
    (A^*)^* = A,
\end{equation}
where $a,b \in \field C$. If $A^*$ exists and coincides with $A$ then $A$ 
is called \emph{symmetric}. Note that if $U: \mathfrak H \to \mathfrak H$ 
is unitary 
then $U^*$ exists and is given by $U^{-1}$. Finally note that $\id^* = \id$. 
Motivated by the familiar case of complex Hilbert spaces one defines 
\begin{equation}
    \begin{array} {rcl}
    \Bounded (\mathfrak H) & := & \{ A \in \End (\mathfrak H) \; | \; 
                                  A^* \textrm { exists } \} \\
    \Unitary (\mathfrak H) & := & \{ U \in \End (\mathfrak H) \; | \;
                                  U \textrm { is unitary } \}.
    \end{array}
\end{equation}
Note that in the particular case where $\mathfrak H$ is indeed a
Hilbert space over the complex numbers then the Hellinger-Toeplitz
theorem ensures that the above definition of $\Bounded (\mathfrak H)$
coincides with the \emph{bounded} operators on $\mathfrak H$, see
e.g.\ \cite [p.~117]{Rud91}. The following lemma is obvious:
\begin{lemma}
Let $\mathfrak H$ be a pre-Hilbert space over $\field C$ then 
$\Bounded (\mathfrak H)$ is a $^*$-algebra with unit element over
$\field C$ and 
$\Unitary (\mathfrak H) \subseteq \Bounded (\mathfrak H)$ is a group. 
Moreover, if $\frac{1}{2} \in \field R$, then any element in 
$\Bounded (\mathfrak H)$ is a $\field C$-linear combination of two
symmetric elements.
\end{lemma}

Now we come to the GNS construction for $^*$-algebras over ordered
rings $\field R$ and the corresponding quadratic extension $\field C$
as above. Let $\mathcal A$ be a $^*$-algebra over $\field C$, i.e.\ an
associative algebra over $\field C$ with an $\field C$-anti-linear
involutive anti-automorphism $^*: \mathcal A \to \mathcal A$. Then a
$\field C$-linear functional $\omega: \mathcal A \to \field C$ is
called \emph{positive} if 
\begin{equation}
\label{PosFunctDef}
    \omega (A^*A) \ge 0
\end{equation}
for all $A \in \mathcal A$. If $\omega$ is positive then the
\emph{Cauchy-Schwarz inequality}
\begin{equation}
    \begin {array} {rcl}
    \omega (A^*B) & = & \cc{\omega (B^*A)} \\
    \omega (A^*B) \cc{\omega (A^*B)} & \le & \omega (A^*A) \omega (B^*B)
    \end {array}
\end{equation}
holds for all $A, B \in \mathcal A$ and implies that the space
\begin{equation}
\label{GelfandIdealDef}
    \mathcal J_\omega := \{A \in \mathcal A \;|\; \omega (A^*A) = 0 \}
\end{equation}
is a left ideal in $\mathcal A$, the so-called 
\emph{Gel'fand ideal}. The quotient space 
$\mathfrak H_\omega := \mathcal A \big/ \mathcal J_\omega$ thus
carries an $\mathcal A$-leftmodule structure given by
\begin {equation}
\label{GNSRepDef}
    \pi_\omega (A) \psi_B := \psi_{AB} ,
\end{equation}
where $\psi_B \in \mathfrak H_\omega$ denotes the equivalence class of
$B$. This representation is called the \emph{GNS representation} of 
$\mathcal A$ on $\mathfrak H_\omega$ induced by $\omega$. Moreover
$\mathfrak H_\omega$ becomes a pre-Hilbert space over $\field C$ by
setting 
\begin{equation}
\label{GNSHermiteanProd}
    \SP{\psi_A, \psi_B} := \omega (A^*B) ,
\end{equation}
which turns out to be a Hermitian product indeed. Finally $\pi_\omega$
is even a $^*$-representation, i.e.\ $(\pi_\omega (A))^*$ always exists
for all $A \in \mathcal A$ and is given by $\pi_\omega (A^*)$. Hence
$\pi_\omega: \mathcal A \to \Bounded (\mathfrak H)$ is a morphism of
$^*$-algebras over $\field C$.

%
%

\section {Formal series and $\lambda$-adic topology}
\label{FormalApp}

In this appendix we shall collect some well-known results on formal
series and the $\lambda$-adic topology. The reader is referred to the
standard algebra textbooks and for Newton-Puiseux and CNP series we
refer to \cite{Ruiz93,BW98a,BNW99a}.

In order to define the formal Laurent, Newton-Puiseux, and completed
Newton-Puisex (CNP) series we first have to specify the allowed
exponents of the formal parameter: Let $S \subset Q$ be a subset with
either a smallest element $q_0 \in S$ or $S = \emptyset$. Then $S$ is
called \emph{CNP-admissable} if $S$ has no accumulation point,
\emph{NP-admissable} if there exists a $N \in \mathbb N$ such that
$N \cdot S \subset \mathbb Z$, and \emph{L-admissable} if 
$S \subset \mathbb Z$, respectively. Now let $V$ be a module over some 
ring $\field R$ and $f: \mathbb Q \to V$ a map. Then one defines the
$\lambda$-support of $f$ by 
$\lsupp f := \{q\in \mathbb Q\;|\; f(q) \ne 0\}$, and the formal
Laurent, Newton-Puiseux, and CNP series with coefficients in $V$ by
\begin{equation}
\label{LNPCNPDef}
    \begin{array} {c}
        \LS V  = \{f: \mathbb Q \to V \; | \; 
                 \lsupp f \textrm{ is L-admissable } \}, \\
        \NP V  = \{f: \mathbb Q \to V \; | \; 
                 \lsupp f \textrm{ is NP-admissable } \}, \\
        \CNP V = \{f: \mathbb Q \to V \; | \; 
                 \lsupp f \textrm{ is CNP-admissable } \},
    \end{array}
\end{equation}
respectively. Observe that 
$V[[\lambda]] \subseteq \LS V \subseteq \NP V \subseteq \CNP V$ are
again $\field R$-modules, namely sub-modules of the 
$\field R$-module of all maps $\mathbb Q \to V$. Elements 
$f \in \CNP V$ are written more familiar as formal series in the
formal parameter $\lambda$
\begin{equation}
\label{CNPseries}
    f = \sum_{q\in \lsupp f} \lambda^q f_q
    \quad
    \textrm{ with }
    \quad
    f_q = f(q).
\end{equation}
The requirement that $\lsupp f$ has in any case a smallest element
if $f \ne 0$ is crucial for the definition of the \emph{order} 
$o(f) := \min (\lsupp f)$ and one sets $o(0) := +\infty$. One
defines the \emph{absolute value} of $f$ by $\varphi (f) := 2^{-o(f)}$ 
and sets $d_\varphi (f, g) := \varphi (f-g)$ for $f, g \in \CNP V$,
which turns out to define an \emph{ultra-metric}, following from the
\emph{strong triangle inequality} $o (f+g) \ge \min (o(f), o(g))$ 
for the order. The induced topology is
called the \emph{$\lambda$-adic topology} and it is well-known that
$V[[\lambda]]$, $\LS V$, and $\CNP V$ are complete with respect 
to this metric, whereas $\NP V$ is dense in $\CNP V$, see e.g.\ 
\cite[Prop.~2]{BW98a}.

The spaces $\field R[[\lambda]]$, $\LS{\field R}$, $\NP{\field R}$,
and $\CNP{\field R}$ have a natural ring structure and $V[[\lambda]]$, 
$\LS V$, $\NP V$, and $\CNP V$ become modules over these rings. In
case when $\field R$ is even a field then $\LS{\field R}$, 
$\NP{\field R}$, and $\CNP{\field R}$ are fields too. We consider now
homomorphisms of such modules. If $\phi_q \in \Hom_{\field R} (V, W)$
is a $\field R$-module homomorphisms from $V$ to $W$ for 
$q \in S \subset \mathbb Q$, where $S$ is an L-, NP-, or
CNP-admissible set, then 
$\sum_{q\in S} \lambda^q \phi_q$ becomes canonically an element of
$\Hom_{\LS{\field R}} (\LS V, \LS W)$, 
$\Hom_{\NP{\field R}} (\NP V, \NP W)$, or
$\Hom_{\CNP{\field R}} (\CNP V, \CNP W)$, respectively, which induces
the following inclusions
\begin{equation}
\label{HomInclusions}
    \begin{array}{c}
        \LS{\Hom_{\field R} (V, W)} 
        \subseteq \Hom_{\LS{\field R}} (\LS V, \LS W) \\
        \NP{\Hom_{\field R} (V, W)} 
        \subseteq \Hom_{\NP{\field R}} (\NP V, \NP W) \\
        \CNP{\Hom_{\field R} (V, W)} 
        \subseteq \Hom_{\CNP{\field R}} (\CNP V, \CNP W). 
    \end{array}
\end{equation}
In the case of formal power series the corresponding inclusion is
known to be indeed an equality, i.e.\ one has 
$\Hom_{\field R} (V, W)[[\lambda]] 
= \Hom_{\field R[[\lambda]]}(V[[\lambda]], W[[\lambda]])$, 
see \cite[Prop.~2.1]{DL88},
but the above three inclusions are in general proper, see e.g.\ 
\cite[App.~A]{BNW99a}. Nevertheless for formal Laurent series one has
the following characterization:
\begin{lemma}
\label{LambdaContLem}
Let $\Phi: \LS V \to \LS W$ be a $\LS{\field R}$-modules
homomorphism. Then $\Phi$ is continuous in the $\lambda$-adic topology 
if and only if $\Phi \in \LS{\Hom_{\field R} (V, W)}$.
\end{lemma}

Let us remember that if $\field R$ is an ordered ring then
$\field R[[\lambda]]$, $\LS{\field R}$, $\NP{\field R}$, and 
$\CNP{\field R}$ become ordered rings in a canonical way: let 
$0 \ne a = \sum_{q \in \lsupp a} \lambda^q a_q \in \CNP{\field R}$
then one defines $a > 0$ if $a_{\min (\lsupp a)} > 0$ in $\field
R$. It is easily verified that $\CNP{\field R}$ is an ordered ring
again and 
$\field R[[\lambda]] \subseteq \LS{\field R}
\subseteq \NP{\field R} \subseteq \CNP{\field R}$ 
are ordered sub-rings. Moreover the topology induced by the order
coincides with the $\lambda$-adic topology, see e.g.\ 
\cite[Prop.~3]{BW98a}.

We shall now briefly remember the definition of the finite topology of
maps and discuss the relation of the $\lambda$-adic topology. Let 
$V, W$ be $\field R$-modules over a ring $\field R$ and consider the
$\field R$-linear morphisms $\Hom_{\field R} (V, W)$. One defines a
topology for $\Hom_{\field R} (V, W)$ by specifying a basis of
neighborhoods of $0 \in \Hom_{\field R} (V, W)$ in the following
way: let
\begin{equation}
    O_{v_1, \ldots, v_k} := \left\{ A \in \Hom_{\field R} (V, W) 
                            \; | \; 
                            A(v_1) = \cdots = A(v_k) = 0 
                            \right\},
\end{equation}
where $k \in \mathbb N$ and $v_1, \ldots, v_k \in V$. This defines a
basis of neighborhoods of $0$ and thus (by translating) a topology on
$\Hom_{\field R} (V, W)$ called the \emph{finite topology}, see e.g.\
\cite{Jac56}. As one
can easy see it coincides with the compact-open topology of maps when
$V$ and $W$ are discretely topologized. Then a sequence 
$(A_n)_{n \in \mathbb N}$ of elements 
$A_n \in \Hom_{\field R} (V, W)$ converges to 
$A \in \Hom_{\field R} (V, W)$ if and only if for all $v \in V$ one
has $A_n v \to Av$ in the discrete topology of $W$ which is the case
if and only if there exists a $N \in \mathbb N$ (depending on $v$)
such that $A_n v = Av$
for all $n \ge N$. Cauchy sequences are defined as usual and clearly 
$\Hom_{\field R} (V, W)$ is complete, i.e.\ any Cauchy sequence
converges. The following example shows that the finite topology is
quite useful and in general it is \emph{strictly coarser} than the
discrete topology:
\begin{example} 
\label{CauchyDiffLocEx}
Let $M$ be a manifold. Then the completion of
the differential operators $\Diff (C^\infty_0 (M))$ on 
$C^\infty_0 (M)$ in the finite topology 
are the local operators $\Loc (C^\infty_0 (M))$.
\end{example}
\begin{proof}
It is straightforward to see that $\Loc (C^\infty_0 (M))$ is
complete in the finite topology and since clearly 
$\Diff (C^\infty_0 (M)) \subseteq \Loc (C^\infty_0 (M))$ we only have
to construct, for a give $L \in \Loc (C^\infty_0 (M))$, a sequence
$D_n$ of differential operators converging to $L$. But this is
essentially Peetre's theorem: let $(O_n, \chi_n)_{n \in \mathbb N}$ be
an approximate identity then by Peetre's theorem 
$D_n := \chi_n L$ \emph{is} a differential operator, since $\chi_n$
has compact support. It follows easily that $D_n \to L$ in the finite
topology since any $f \in C^\infty_0 (M)$ has support in some $O_n$.
\end{proof} 
\begin{remark}
Since the definition of the finite topology as well as the definition
of differential operators on an associative, commutative algebra is
purely algebraic we note that this example provides a method to
\emph{define} local operators in general, namely as completion
of the differential operators in the finite topology. Of course the
same can be done for multidifferential operators.
\end{remark}

Let us now investigate the connection between the finite topology of
$\Hom_{\field R} (V, W)$ and the `strong operator topology' of
$\LS{\Hom_{\field R} (V, W)}$ which is defined by the following basis
of neighborhoods of $0$: let $\epsilon > 0$ and 
$v_1, \ldots, v_k \in V$ then we define
\begin{equation}
    O_{v_1, \ldots, v_k: \epsilon} :=
    \left\{ A \in \LS{\Hom_{\field R} (V, W)} \; | \;
    \forall l=1, \ldots, k: \; \varphi ( A (v_l)) < \epsilon \right\}
\end{equation}
which clearly determines a topology. Here we used the $\lambda$-adic 
absolute value $\varphi$ on $W$. Clearly a sequence 
$A_n \in \LS{\Hom_{\field R} (V, W)}$ converges to some 
$A \in \LS{\Hom_{\field R} (V, W)}$ in the strong operator topology
if and only if for all $v \in \LS V$ the sequence $A_n v$ converges
$\lambda$-adically to $Av$ which motivates the name of this topology.
Note that the $\lambda$-adic topology of 
$\LS{\Hom_{\field R} (V, W)}$ is (in general strictly) finer that the
strong operator topology, see Example~\ref{DiffLocEx}. 
The following proposition shows that the
finite topology of $\Hom_{\field R} (V, W)$ and the strong operator
topology of $\LS{\Hom_{\field R} (V, W)}$ fit together very
naturally:
\begin{proposition}
\label{CauchyLocProp}
Let $V, W$ be $\field R$-modules and let 
$\Diff \subseteq \Loc \subseteq \Hom_{\field R} (V, W)$ be subspaces
such that $\Loc$ is the completion of $\Diff$ in the finite
topology. Then $\LS\Loc$ is the completion of $\LS\Diff$ in the
strong operator topology of $\LS{\Hom_{\field R}(V, W)}$.
\end{proposition}
\begin{proof}
We shall first show the following lemma which is a particular case of
the proposition:
\begin{lemma}
$\LS{\Hom_{\field R} (V, W)}$ is complete in the strong operator
topology.
\end{lemma}
\begin{innerproof}
Let $A_n \in \LS{\Hom_{\field R} (V, W)}$ be a Cauchy sequence. Since
$\LS W$ is complete in the $\lambda$-adic topology we observe that for
any $v \in \LS V$ the sequence $A_n v$ is Cauchy and thus convergent
in $\LS W$. Thus $Av := \lim_n A_n v$ clearly defines a 
$\LS{\field R}$-linear map $A: \LS V \to \LS W$. It thus remains to
show that $A$ is $\lambda$-adically continuous and that $A_n \to A$ in
the strong operator topology.
To this end we write each $A_n$ as 
$A_n = \sum_{r=o_n}^\infty \lambda^r A^{(r)}_n$ where 
$o_n = o (A_n) \in \mathbb Z$ is the order of $A_n$. We now assume
that the orders $o_n$ are \emph{not} uniformly bounded from below and
lead this assumption to a contradiction. We may thus even assume that
$o_n = -2n$ by choosing a suitable sub-sequence (the factor $2$ is
only for technical reasons and we also could have assumed 
$o_n = -2n-1$). Hence $A^{(-2n)}_n \ne 0$. Choose $v_0 \in V$ such that 
$A^{(0)}_0 v_0 \ne 0$. Then we can find $v_1 \in V$ such that 
$A_1 (v_0 + \lambda v_1)$ has a non-vanishing term in order 
$\lambda^{-1}$: either $A^{(-1)}_1 v_0 \ne 0$ then choose $v_1 = 0$ or
$A^{(-1)}_1 v_0 = 0$ then choose $v_1$ such that 
$A^{(-2)}_1 v_1 \ne 0$ which is possible by assumption. Inductively we
can find $v_0, v_1, \ldots, v_n, \ldots \in V$ such that 
$A_n (v_0 + \lambda v_1 + \cdots + \lambda^n v_n)$ has a non-vanishing  
term in order $\lambda^{-n}$. Defining now 
$v := \sum_{n=0}^\infty \lambda^n v_n$ we observe that due to the
$\lambda$-adic continuity $A_n v$ has a non-vanishing term of order
$\lambda^{-n}$ too, whence $-2n \le o (A_n v) \le -n$. But this is in
contradiction to $A_n v \to Av$. Thus we conclude that the orders
$o_n$ are bounded from below by some $N \in \mathbb Z$. It clearly
follows that $o (Av) \ge N + o(v)$ for all $v \in \LS V$,
whence $A$ is $\lambda$-adically
continuous. Then $A_n \to A$ in the strong operator topology follows
by construction of $A$.  
\end{innerproof}

\noindent
\textsc{End of the proof of the proposition:}
Now let $L_n \in \LS\Loc$ be a Cauchy sequence with respect to the
strong operator topology. Then $L_n \to L$  with
some $L \in \LS{\Hom_{\field R} (V, W)}$ by the lemma.
Moreover, we know due to the last lemma that the orders of $L_n$ and
$L$ are bounded from below by some $N \in \mathbb Z$. Considering 
$v \in V$ we have $L_n v \to Lv$ in the $\lambda$-adic topology. In
lowest order $N$ this implies that the order $N$ of $L_n$ converge
to the order $N$ of $L$ in the finite topology of 
$\Hom_{\field R} (V, W)$ whence by Cauchy completeness of $\Loc$ we
find that the lowest order of $L$ is in $\Loc$. Now an easy induction
shows that indeed $L \in \LS\Loc$ proving the completeness of $\LS\Loc$
in the strong operator topology. Consider finally $L \in \LS\Loc$
written as $L = \sum_{r=N}^\infty \lambda^r L^{(r)}$ with 
$N \in \mathbb Z$ and let $D^{(r)}_n \in \Diff$ be a sequence such
that in the finite topology we have $D^{(r)}_n \to L^{(r)}$. We claim
that $D_n := \sum_{r=N}^\infty \lambda^r D^{(r)}_n \to L$ in the
strong operator topology which is indeed the case as an easy
verification shows. Thus the proposition is shown.
\end{proof}

%
%

\section*{Acknowledgments}

I would like to thank Martin Bordemann for many useful discussions
and suggestions, in particular concerning the finite
topologies. Moreover, I would like to thank Nikolai Neumaier and Alan 
Weinstein for valuable discussions. Finally, warm hospitality of
the Math Department of UC Berkeley, where this work has been finished,
and financial support of the Studienstiftung des deutschen Volkes is
gratefully acknowledged.

%
%

\end{document}